\documentclass{amsart}

\usepackage{amscd}
\usepackage{verbatim}
\usepackage{amssymb}

\theoremstyle{plain}
\newtheorem{Thm}{Theorem}

\newtheorem{Lem}{Lemma}
\newtheorem{Prop}{Proposition}

\theoremstyle{definition}
\newtheorem{Def}{Definition}

\theoremstyle{remark}

\numberwithin{equation}{section}
\newcommand{\PS}[1]{{\Bbb P} ^{#1}}                  
\newcommand{\Op}[1]{{\mathcal O}_{{\Bbb P}^{#1}}} 
\newcommand{\sexs}[3]{0 \longrightarrow #1 \longrightarrow #2
\longrightarrow #3 \longrightarrow 0}
\newcommand{\res}[1]{\!\!\mid_{#1}}
\newcommand{\mL}{{\mathcal L}}
\newcommand{\mO}{{\mathcal O}}
\newcommand{\znums}{{\mathbf Z}}
\newcommand{\cnums}{{\mathbf C}}

\newcommand{\sheafext}{\operatorname{\mathcal{E}\!\mathit{xt}}}
\newcommand{\sheafhom}{\operatorname{\mathcal{H}\!\mathit{om}}}
\newcommand{\sheafend}{\operatorname{\mathcal{E}\!\mathit{nd}}}
\newcommand{\sheaftor}{\operatorname{\mathcal{T}\!\mathit{or}}}
\newcommand{\mP}{{\mathbf P}}
\newcommand{\OS}{{\mathcal O}_{S}}
\newcommand{\mM}{{\mathcal M}}
\newcommand{\mY}{{\mathcal Y}}
\newcommand{\mI}{{\mathcal I}}
\newcommand{\mK}{{\mathcal K}}
\newcommand{\mF}{{\mathcal F}}
\newcommand{\mE}{{\mathcal E}}
\newcommand{\mS}{{\mathcal S}}
\newcommand{\mW}{{\mathcal W}}
\newcommand{\mQ}{{\mathcal Q}}
\newcommand{\mmP}{{\mathcal P}}
\newcommand{\Sym}{\operatorname{Sym}}
\newcommand{\Ker}{\operatorname{Ker}}

\begin{document}

\title[Stable Bundles]{On Stable Bundles of Ranks 2 and 3\\
   on $ {\Bbb P} ^{3}$}
\author{Al Vitter}
\address{Department of Mathematics\\
   Tulane University \\
   New Orleans, La. 70118\\
   USA}
\email{vitter@math.tulane.edu}
\thanks{}

\date{September 25, 2003}


\maketitle
\begin{abstract}
We study rank 3 stable bundles E on $\PS{3}$ as extensions of a
line bundle $\mL$ on a smooth surface $S\subset \PS{3}$ by
$\overset{3}{\oplus}\Op{3}(-\nu).$ In most cases, S (the
dependency locus of three sections of $E(\nu)$) lies in the
Noether-Lefschetz locus. We give a detailed analysis when S
contains a line L and $\mL$ is constructed from divisors of the
form $aL+bC$ for $H=L+C$ a hyperplane section of S. We study
the parameter space of this construction and compare it to the
full (Gieseker-Maruyama)moduli space. We also analyze the
situation when $\mL$ is a power of the hyperplane bundle.

The same approach is used to study rank 2 stable bundles on $\PS{3}.$
\end{abstract}

\section{Introduction} \label{S:intro}

The purpose of this paper is to begin a study of stable vector
bundles of rank three on three dimensional projective space. Our
approach is to express such a bundle E (normalized so that $c_{1}=0,-1$ or -2)
 as an extension

\begin{equation} \label{E:seq1}
\begin{CD}
0 @>>> \overset{3}{\oplus}\Op{3}(-\nu) @>> \sigma > E
@>>> {j_{S}}_{*}\mL @>>> 0
\end{CD}
\end{equation}

\noindent  for $S \subset \PS{3}$ a smooth surface of degree $k=3\nu +c_{1}$ and $\mL$
a line bundle on S, using Serre's Theorem A and the Kleiman Transversality
Theorem. We study E through S, $\mL$, and the extension class $\tau$
of (\ref{E:seq1}) which appears in the dual sequence

\begin{equation} \label{E:seq2}
\begin{CD}
0 @>>> E^{*} @>>\sigma^{t} >
\overset{3}{\oplus}\Op{3}(\nu) @>>\tau >
{j_{S}}_{*}{\mathcal L}^{*}(k) @>>> 0. \\
\end{CD}
\end{equation}

\noindent Chern class calculations show that, in most cases, S must
belong to the Noether-Lefschetz locus, that is, it must support a
line bundle not equal to a power of the hyperplane bundle.

To produce examples, we reverse the above procedure and start with
$c_{1}\in \{0,-1,-2\}$,
$\nu \in \znums_{+}$, a surface $S\subset\PS{3}$ of degree $k=3\nu+c_{1}$,
and a line bundle $\mL$ on S and
consider extensions (\ref{E:seq1}). We make a
detailed study of the case where the surface contains a line L and the line bundles
are constructed from divisors of the form $aL+bC$ for $L+C$ a
hyperplane section of S containing L, C a curve of degree k-1, and $a,b\in \znums$
(Section~\ref{S:examples of rank three, part II}).
It is determined
when the resulting coherent sheaf E is locally free and (modulo one unresolved case)
 when it is stable (Theorem~\ref{Thm:main theorem rank 3}). We
count the moduli of our construction (Proposition~\ref{Prop:parameters rank 3})
by proving that the correspondence $(E,\sigma)\leftrightarrow (S,L,\tau)$
is 1-to-1. Then we estimate the
dimension of the component of the full moduli space containing E,
$\mM$ (Theorem~\ref{Thm:moduli for k large}).
When the degree of S is 2 or 3, $dim\mM$ is determined exactly
and we can conclude, in many cases, that our examples form a subset $\mY$ of $\mM$ of
equal dimension
and that $\mM$ is smooth at E
(Theorem~\ref{Thm:moduli degS=2} and Theorem~\ref{Thm:moduli
degS=3}). For arbitrary k, we give a separate analysis of the special case
where the line bundle $\mL$ is a power of the hyperplane bundle
(Section~\ref{S:examples of rank 3}) and show that the corresponding
space of parameters $\mY$ is an open subscheme of $\mM.$ We address the general problem
of putting a scheme structure on the parameter space $\mY$ in
Section~\ref{S:scheme}.

The examples we construct and study provide  evidence for the
general problem of determining the dimension of the moduli space
of stable bundles when the base variety has dimension $\geq 3.$
For $E$ a rank $r$ stable bundle on  a smooth projective variety $X$
and $\mM$ the corresponding moduli space (see
Section~\ref{S:pre}), $T\mM_{E} \cong H^{1}(X;{\sheafend}_{0}E)$ (${\sheafend}_{0}E$
is the bundle of trace-free endomorphisms of $E.$) and

\begin{equation} \label{E:intro-eq3}
h^{1}(X;{\sheafend}_{0}E)-h^{2}(X;{\sheafend}_{0}E) \leq \text{
dim}_{E}\mM \leq h^{1}(X;{\sheafend}_{0}E).
\end{equation}

\noindent The expected dimension of $\mM$ is defined by

\begin{equation} \label{E:intro-eq5}
ed(\mM)\equiv h^{1}(X;{\sheafend}_{0}E)-h^{2}(X;{\sheafend}_{0}E).
\end{equation}

\noindent When $X$ is a surface, Riemann-Roch calculates

\begin{equation} \label{E:intro-eq4}
ed(\mM)\equiv 2rc_{2}(E)-(r-1)c_{1}(E)^{2}-(r^{2}-1)\chi(\mO_{X}).
\end{equation}

\noindent Also for the surface case, important work by Gieseker
and
Li (\cite{Gieseker-Li-moduli-2} and \cite{Gieseker-Li-moduli}), and
O'Grady \cite{O'Grady-moduli} implies that, for $c_{2}(E)$ large enough (with $c_{1}(E)$
fixed), $\mM$ is irreducible, generically smooth, and of dimension
$ed(\mM)$, and, on a Zariski open subset of $\mM$, $h^{2}(X;{\sheafend}_{0}E)=0.$

When the base variety has dimension $\geq 3$, no results of this
type have been proven. And there is no expression for $ed(\mM)$ in
terms of chern classes (like (\ref{E:intro-eq4}))- because of the
higher dimensional groups $h^{i}(X;{\sheafend}_{0}E)$, $i\geq 3.$
By varying the discrete parameters in the examples of
Section~\ref{S:examples of rank 2}, Section~\ref{S:examples of
rank 3}, and Section~\ref{S:examples of rank three, part II}, one
finds many bundles E for which $\text{dim }\mM$ is much larger than $ed(\mM)$
and for arbitrarily large $c_{2}(E).$ For these examples, $h^{2}(X;{\sheafend}_{0}E)$
is in fact much larger than $ed(\mM).$ One could ask whether the term "expected
dimension" should be applied to (\ref{E:intro-eq5}) when the base
manifold has dimension three or greater. The problem
remains:Understand $\text{dim }\mM$ for stable bundles over
smooth varieties of dimension $\geq 3.$

The technical backbone of this paper's theorems consists of the
intersection properties of L and C on S and results on the
cohomology of the line bundles $\OS(iL+jC)$ (Section~\ref{S:surfaces with lines}).

The same methods are also applied to stable rank two bundles on $\PS{3}$
(Section~\ref{S:sb,r2} and Section~\ref{S:examples of rank 2}). In
general, the examples produced from surfaces containing a line
seem to comprise a higher codimension subset of $\mM$ than in
the rank three case.

Our approach can also be used to discuss stable bundles $E\rightarrow X$
of various ranks on other smooth projective varieties X. This will
be the subject of future papers.

It is a pleasure to thank Jim Bryan and Bob Friedman for helpful
conversations.

\bigskip

\section{Preliminaries}
\label{S:pre}
By a stable bundle we shall mean Mumford-stable (or $\mu-$
stable), that is

\begin{Def}
Let X be a smooth projective variety of dimension n, $\mO_{X}(1)$
a very ample line bundle on X,
and H a corresponding hyperplane section of X. A coherent torsion-free rank r sheaf E
on X is called stable ( resp. semistable) if, for any subsheaf $F\subset E$ of rank
$r_{1}<r$, $r_{1}^{-1} c_{1}(F)\cdot H^{n-1}<r^{-1}c_{1}(E)\cdot H^{n-1}$
\quad (resp. $\leq$ ).
\end{Def}

\begin{Def}
A coherent torsion-free rank r sheaf E on X is called Gieseker-stable
(resp. Gieseker-semistable) if, for any proper subsheaf $F\subset E$
of rank $r_{1}$, $r_{1}^{-1}\chi (X;F(l))<r^{-1}\chi (X;E(l))$
(resp. $\leq$ ) for $l\gg 0$, where $\chi (X;E(l))$ is the
Hilbert polynomial of E.
\end{Def}

There is a coarse moduli space (\cite[page 153]{Friedman} and \cite[page 38 and
chapter 4]{Huybrechts-Lehn}) for the Gieseker-semistable sheaves on X
with fixed Hilbert polynomial, a projective scheme $\mM$
whose closed points correspond to the S-equivalence classes
(\cite[page 22]{Huybrechts-Lehn})of
Gieseker-semistable sheaves on X. From the definitions and
Riemann-Roch it follows that stable $\Rightarrow$ Gieseker-stable $\Rightarrow$
Gieseker-semistable $\Rightarrow$ semistable. The stable sheaves
with fixed $\chi (X;E(l))$ form an open subset of
$\mM.$

The Riemann-Roch formula \cite[Append. A, sec.4]{Hartshorne}
for a rank r coherent sheaf E on $\PS{3}$ is

$$
\chi(\PS{3};E)=r+\frac{11}{6}c_{1}+(c_{1}^{2}-c_{2})+
\frac{1}{6}(c_{1}^{3}-3c_{1}c_{2}+3c_{3}).
$$

\noindent Here we have identified the chern classes
$c_{i}(E)=c_{i}$ with integers using the positive generator $\omega_{0}$
of $H^{2}(\PS{3};\znums)$ and the generators $\omega_{0}^{i}$ of
$H^{2i}(\PS{3};\znums)$ i=0 to 3.

The Grothendieck-Riemann-Roch formula for a closed embedding of
smooth varieties $f:X\rightarrow Y$ and a coherent sheaf E on X
(\cite[Chapter 15]{Fulton} and \cite[Append. A, sec.4]{Hartshorne})is

$$
ch (f_{*}E) = f_{*}[ch(E){td(N_{X\mid Y})}^{-1}].
$$

We make frequent use of the Kleiman Transversality Theorem
(\cite{Kleiman} and \cite[Thm10.8]{Hartshorne}): Let X be a homogeneous variety with group
variety G over an algebraically closed field k of characteristic
0. Let $f:Y\rightarrow X$ and $\phi:Z\rightarrow X$ be morphisms of
nonsingular varieties Y,Z to X. For any $g \in G(k)$, let
$Y^{g}$ be Y with the morphism $g \circ f$ to X. Then
there is a nonempty (Zariski)open subset $U\subset G$ such that
for every $g \in U(k)$, $Y^{g}\times _{X} Z$ is
nonsingular and either empty or of dimension exactly $dimY+dimZ-dimX.$

We use Kleiman Transversality in the following situation (see \cite[page 121]{Huybrechts-Lehn}).
 Let Y be a smooth
projective variety of dimension n and E a rank r vector bundle on Y which is
globally generated. Set $H\equiv H^{0}(Y;E).$
For X the grassmannian of r-dimensional quotient spaces of H,
evaluation of sections defines a regular map $f:Y\rightarrow X$
such that, $E\cong f^{*}\mQ$ for $\mQ$ the tautological quotient
bundle on X. For any k sections of E, $\sigma_{j}$ j=1 to k, which
generate a k-dimensional subspace V of H, and $0\leq l\leq k$, set
$Y_{l}\equiv \{y\in Y \mid \text{ dim } span\{\sigma_{1}(y),\dots,\sigma_{k}(y)\}
\leq l\}.$
Define $Z_{l}\equiv \{ H/K \in X \mid \text{ dim } K\cap V \geq k-l \}$,
a Schubert variety of codimension $(k-l)(r-l).$  $Z_{l}$ is smooth away from
$Z_{l-1}$ and $Y_{l}= f^{-1}Z_{l}.$ The group G is $GL(dim H,\cnums)$ acting on X.
Now Kleiman Transversality implies that, for
generic $\sigma_{j}$ j=1 to k, $Y_{l}$ is of
codimension $(k-l)(r-l)$, empty if $(k-l)(r-l)>n$, and the singular locus of $Y_{l}$ is
of codimension $(k-l-1)(r-l-1)$, empty if $(k-l-1)(r-l-1)>n.$

\bigskip

\section{Stable Bundles of Rank 2 on $\PS{3}$}
\label{S:sb,r2}

Let $ E \longrightarrow {\Bbb P}^{3} $ be a rank 2 normalized
bundle ($c_{1}=0  \text{ or }  -1$). For $\nu$ large enough, $E(\nu)$
is globally generated and a generic $\sigma =
(\sigma_{1},\sigma_{2})\in \overset {2}{\oplus} H^{0}({\Bbb
 P} ^{3};E(\nu)) $ gives an exact sequence

\begin{equation} \label{E:sequence for E}
\begin{CD}
o @>>> \overset{2}{\oplus}\Op{3}(-\nu) @>> \sigma_{1}\oplus \sigma_{2}> E
@>>> {j_{S}}_{*}\mL @>>> 0
\end{CD}
\end{equation}

\noindent (For fixed chern classes $ \exists \nu_{0} \in {\Bbb Z}_{+} $ so
that this holds $ \forall \nu \geq \nu_{0} $ and all semistable E,
since this family is bounded \cite[Thm 3.3.7]{Huybrechts-Lehn}
\cite[Thm 1.1]{Simp1}). By Kleiman transversality,
the generic $ \sigma $ produces a degeneracy locus $ S =
Z_{\sigma_{1}\wedge \sigma_{2}}$ which is a smooth hypersurface $ S
\mathop{\hookrightarrow}\limits_{j_{S}}{\Bbb P}_{3} $ of degree $
k=2\nu+c_{1} $ , a line bundle ${\mL}$ on S, and zero sets $
Z_{\sigma_{j}}$ j=1,2 which are smooth curves of degree $
c_{2}(E(2)) = c_{2}(E)+c_{1}\nu+\nu^{2}$. It follows that, though
the $ Z_{\sigma_{j}} $ need not be irreducible, their components
are mutually disjoint. This gives a basepoint-free pencil of
curves on S, $Z_{t_{1}\sigma_{1}+t_{2}\sigma_{2}} \quad
[t_{1},t_{2}]\in {\Bbb {P}}^{1} $ and thus a regular map $S
\longrightarrow {\Bbb {P}}^{1}$. Therefore

\bigskip
\begin{it}
\noindent S belongs to the Noether-Lefschetz locus, i.e. it supports a line
bundle not equal to a power of the hyperplane bundle.
\end{it}

\bigskip
Applying $Hom_{\Op{3}}(\quad ,\Op{3})$ to (\ref{E:sequence for E})
gives

\begin{equation} \label{E:dual equation for E}
\begin{CD}
0 @>>> E^{*} @>>(\sigma_{1}^{t},\sigma_{2}^{t})>
\overset{2}{\oplus}\Op{3}(\nu) @>>\tau_{1}\oplus\tau_{2}>
{j_{S}}_{*}{\mathcal L}^{*}(k) @>>> 0 \\
\end{CD}
\end{equation}

\noindent for $\tau=(\tau_{1},\tau_{2}) \in
\overset{2}{\oplus} H^{0}(S;{\mathcal L}^{*}(\nu+c_{1}))$. I explain
why, after possibly multiplying $\tau_{1}$ and $\tau_{2}$ by the
same non-zero constant,

\begin{equation} \label{E:tau}
\tau_{1}=\sigma_{2}\res{S}\wedge \qquad
\tau_{2}=-\sigma_{1}\res{S}\wedge.
\end{equation}

\noindent View ${\mathcal L}$ as the quotient sheaf $E \diagup
im(\sigma_{1}\oplus\sigma_{2})$. For $g=(g_{1},g_{2})\in
\overset{2}{\oplus}\Op{3}(\nu)$ ,set
$T(g)\equiv(g_{1}\sigma_{2}-g_{2}\sigma_{1})\res{S}\wedge$ and
apply to $[f]\in \mL$, $f\in E$ to get \linebreak
$(g_{1}\sigma_{2}-g_{2}\sigma_{1})\wedge f \res{S} \in
{\mathcal O}_{S}(k).$ Note that this is well-defined independent of
$f\in[f]$ and applied to $g=(\sigma_{1}^{t},\sigma_{2}^{t})(\psi)=
(\psi(\sigma_{1}),\psi(\sigma_{2}))$ for $\psi\in E^{*}$ gives
\linebreak
$(\psi(\sigma_{1})\sigma_{2}-\psi(\sigma_{2})\sigma_{1})\res{S}\wedge =
\iota_{\psi}(\sigma_{1}\wedge\sigma_{2})\res{S}\wedge =0.$ It
follows that $T=\tau_{1}\oplus\tau_{2}$ up to non-zero constant
multiple.

 $Z_{\sigma_{1}}=Z_{\tau_{2}}, \quad Z_{\sigma_{2}}=Z_{\tau_{1}},
 \quad Z_{\sigma_{1}}\cdot Z_{\sigma_{2}}=0$, and so
 $Z_{\tau_{1}}\cdot Z_{\tau_{2}} =0$ (which also follows from
 (\ref{E:dual equation for E})).
 Therefore $c_{1}({\mathcal L}^{*}(\nu+c_{1})
 =(\nu+c_{1})\omega_{0}-c_{1}({\mathcal L})$ gives

\begin{equation} \label{E:chern1}
((\nu+c_{1})\omega_{0}-c_{1}({\mathcal L}))^{2}=0 \quad
\text{(intersection on S) i.e.}
\end{equation}

\begin{equation} \label{E:chern2}
(\nu+c_{1})^{2}(2\nu+c_{1})-2(\nu+c_{1})\omega_{0} \cdot c_{1}(\mL) +
c_{1}(\mL)^{2}=0.
\end{equation}

Now assume that E is stable. If $\mL$ has the form $\mL={\mathcal O}_{S}(l)$,
$h^{0}(\PS{3};E)=0$ implies $l<0$ but then (\ref{E:chern1}) gives $(\nu+c_{1}-l)^{2}=0$
which is impossible. Therefore

\bigskip
\begin{it}
\noindent For E stable, $\mL \ne \mO _{S}(l)$ for any $l.$
\end{it}

\bigskip
Applying Grothendieck Riemann-Roch to ${j_{S}}_{*} \mL$ \cite{Fulton} gives

\begin{equation}   \label{E:second chern class of E}
c_{2}(E)=(\nu+c_{1})^{2} - \omega_{0} \cdot c_{1}(\mL) \quad
\text{(intersection on S)}
\end{equation}

\noindent and reproves (\ref{E:chern2}) from the fact that $c_{3}(E)=0.$

To construct some concrete bundles E, reverse the above procedure,
begin with a given $\nu \in \znums_{+} $ , a smooth $S \hookrightarrow \PS{3}$
of degree $k=2\nu + c_{1}\quad (c_{1} =0 \text{ or } -1)$ , and a
line bundle $\mL$ on S and consider extensions

\begin{equation} \label{E:extension}
\sexs{\overset{2}{\oplus} \Op{3}(-\nu)}{E}{{j_{S}}_{*}\mL}.
\end{equation}

These are classified by $Ext^{1}(\PS{3};{j_{S}}_{*}\mL,\overset{2}{\oplus}\Op{3}(-\nu))
\cong \overset{2}{\oplus}H^{0}(S;\mL^{*}(\nu+c_{1}))$;
we want to determine which extensions are locally free. Applying
${\sheafhom}_{\Op{3}}(\quad ,\Op{3})$ to (\ref{E:extension}) gives

\begin{equation} \label{E:dual extension}
\begin{CD}
0 @>>> E^{*} @>>> \overset{2}{\oplus}\Op{3}(\nu)
@>>\tau_{1}\oplus\tau_{2}> \mL^{*}(k) @>>>
{{\sheafext}^{1}}_{\Op{3}}(E,\Op{3}) @>>> 0.
\end{CD}
\end{equation}

\noindent E is locally free iff ${{\sheafext}^{1}}_{\Op{3}}(E,\Op{3})=0$
iff $\mL^{*}(\nu+c_{1})$ is globally generated by $\tau=(\tau_{1},\tau_{2})$
which is the extension class mentioned above. It follows that

\bigskip
\begin{it}
\noindent The generic extension (\ref{E:extension}) is locally free iff $\mL^{*}(\nu+c_{1})$
is globally generated (necessarily by two sections). In this
case, $((\nu+c_{1}) \omega_{0}-c_{1}(\mL))^{2}=0.$
\end{it}

\bigskip
Recall that a rank 2 bundle E on $\PS{3}$ is stable iff
$h^{0}(\Op{3};E)=0$ and semistable $(c_{1}=0 \text{ case})$ iff $h^{0}(\Op{3};E(-1))=0$
 \cite[pages 165--166]{OSS}. This gives

\bigskip
\begin{it}
\noindent E of the form (\ref{E:extension}) is stable iff $\nu >0$ and $h^{0}(S;\mL)=0.$\\
E is semistable ($c_{1}=0$ case) iff $\nu\geq 0$ and $h^{0}(S;\mL(-1))=0.$
\end{it}

\bigskip
Now S must be chosen from the Noether-Lefschetz locus. The
hypersurfaces of $\PS{3}$ of degree k are parametrized by a $\PS{N_{k}}$
for $N_{k}=\binom{k+3}{3}-1.$ M. Noether stated and Lefschetz
proved that there is a countable union of subvarieties $NL\subset\PS{N_{k}}$
such that $S\notin NL$ implies $Pic(S)\cong \znums$ is generated
by $\mO_{S}(1)$. See \cite{GH-NL}
for a modern proof and also \cite{CHM} and \cite{Lopez} for
interesting properties, references, and questions about the
Noether-Lefschetz locus. The component of NL of smallest
codimension k-3 (and the only such component) consists of the
surfaces in $\PS{3}$ containing a line \cite{Green-NL1,Green-NL2,Voisin-NL1}.

\section{Surfaces in $\PS{3}$ Containing a Line}
\label{S:surfaces with lines}

Let $S \in \PS{3}$ be a smooth degree k surface ($k \geq
2$)containing a line L. Denote the pencil of hyperplane sections
of S containing L by $H_{t} \quad t \in \PS{1}$ and let H be a general
hyperplane section (not containing L).

\begin{Lem} \label{L:intersection}
For $H_{t} =L+C_{t}$, degree $C_{t} = k-1$ and, using intersection on S,
\begin{align}
L^{2} &= 2-k \notag \\
C_{t} \cdot L &= k-1 \notag \\
C^{2}_{t} &= 0. \notag
\end{align}
Furthermore, the generic $C_{t}$ is irreducible and smooth and the pencil $\{ C_{t} \}$
is base point free and thus gives a regular map
$S\overset{\pi}{\longrightarrow}\PS{1}$.
\end{Lem}

\begin{proof}
The genus formula applied to L gives
$0=1+\frac{1}{2}(L^{2}+K_{S}\cdot L) =1+\frac{1}{2}(L^{2}+k-4)$
i.e. $L^{2}=2-k$.

\begin{align}
H^{2}_{t} &= k = L^{2}+2L\cdot C_{t}+C^{2}_{t} \notag \\
2k-2 &= 2L\cdot C_{t}+C^{2}_{t} \notag \\
H_{t}\cdot C_{t} &= k-1 = L\cdot C_{t}+C^{2}_{t}
\end{align}

\noindent Subtraction gives $L\cdot C_{t}=k-1$ and so
$C^{2}_{t}=0$. The base locus of $\{C_{t}\}$ is contained in L and
therefore is finite. Now $0=C^{2}=C_{t_{1}}\cdot C_{t_{2}}\geq 0$
implies that $\{C_{t}\}$ is base point free.
Use C to denote an arbitrary $C_{t}$ and consider

\begin{equation} \label{E:line bundle of C}
\sexs{\mO_{S}}{\mO_{S}(C)}{\mO_{C}(C)}.
\end{equation}

\noindent Since $H^{1}(S;\mO_{S}(j))=0 \quad \forall j$ (as follows from the
cohomology sequence of $\sexs{\Op{3}(j-k)}{\Op{3}(j)}{\mO_{S}(j)}$),
and $\mO_{C}(C)=\mO_{C}$, the cohomology sequence of
(\ref{E:line bundle of C}) gives

$$
\sexs{H^{0}(S;\mO_{S})}{H^{0}(S;\mO_{S}(C))}{H^{0}(C;\mO_{C})}
$$

\noindent which shows that $h^{0}(S;\mO_{S}(C))=2$ and that $\mO_{S}(C)$ is
globally generated.
Bertini's theorem implies that the generic $C_{t}$ is smooth.
If $C_{t}$ had two distinct irreducible components, they must be
disjoint by smoothness. But this is impossible because they are
both contained in the same plane. Therefore $C_{t}$ is irreducible.
\end{proof}

Let $\mO_{S}(aL+bC)(j)\equiv \mO_{S}(aL+bC)\otimes \mO_{S}(j)$ for
a, b, j $\in \znums$ and note that this is isomorphic to $\OS((a-b)L)(b+j)$
and to $\OS((b-a)C)(a+j).$ We will make frequent use of these
isomorphisms and the

\begin{Lem}  \label{Lem:Main Lemma}
For a, b, j $\geq 0$
\begin{align}
i)\quad &H^{0}(S;\OS(-aL)(j))=0 \text{ iff }a>j. \notag \\
ii)\quad &H^{0}(S;\OS(-bC)(j))=0 \text{ iff } b>j. \notag \\
iii) \quad &h^{0}(S;\OS(bC))=b+1 \text{ and } \OS(bC)\text{ is globally
generated.} \notag \\
iv) \quad &H^{1}(S;\OS(-aL)(-j))=0 \text{ iff } j>(a-1)(k-2)
\text{ or } j=0, a=1 \text{ or } a=0.
\notag \\
v) \quad &H^{1}(S;\OS(-bC)(-j))=0 \text{ iff } j>0, \text{ or } j=0,
b=1 \text{ or } b=0. \notag \\
vi) \quad &\text{For } j>k-4 \text{ and } b \geq 0, h^{0}(S;\OS(bC)(j))= \binom{j+3}{3}
-\binom{j-k+3}{3} \notag \\
\quad &+b[\binom{j+2}{2}-\binom{j-k+3}{2}]. \notag
\end{align}
\end{Lem}

\begin{proof}
For $\sigma \in H^{0}(S;\OS(-aL)(-j))$ ,choose C which is not an
irreducible component of $Z_{\sigma}.$ Then $0\leq Z_{\sigma} \cdot C = (-aL+jH) \cdot
C=-a(k-1)+j(k-1)$ so $a \leq j.$ If $a\leq j, \OS(-aL)(j)$ clearly
has global sections so i) holds.

To prove ii) first note that $\OS(-bC)(j)$ also clearly has
sections if $b \leq j.$ If $b>j$ suppose $\sigma \in
H^{0}(S;\OS(-bC)(j)).$ If $L\nsubseteq Z_{\sigma}, 0\leq L \cdot Z_{\sigma} =L \cdot
(-bC+jH)=-b(k-1)+j
<0$ , a contradiction. Therefore $L \subset Z_{\sigma}$ and so
$\OS(-bC-L)(j)= \OS(-(b-1)C)(j-1)$ has a global section. Repeating
this argument gives that $\OS(-(b-j)C)$ has a non-zero global
section, which is not true.

Note that $\OS(bC)\cong \OS(C)^{\otimes b}$ is globally generated
because $\OS(C)$ is. This and the cohomology sequence of

$$
\sexs{\OS((j-1)C)}{\OS(jC)}{\mO_{C}(jC) \cong \mO_{C}}
$$

\noindent gives

$$
\sexs{H^{0}(S;\OS((j-1)C))}{H^{0}(S;\OS(jC))}{H^{0}(C;\mO_{C})}
$$

\noindent and iii) follows by induction.

The group $H^{1}(S;\OS(-j))$ vanishes for all j. Careful examination
of the cohomology sequences

\begin{multline}
0 \longrightarrow H^{0}(S;\OS(-iL)(-j)) \longrightarrow H^{0}(S;\OS(-(i-1)L)(-j)) \\
\longrightarrow H^{0}(L;\mO_{L}((k-2)(i-1)-j)) \longrightarrow
H^{1}(S;\OS(-iL)(-j))\\
\longrightarrow H^{1}(S;\OS(-(i-1)L)(-j))
\longrightarrow H^{1}(L;\mO_{L}((k-2)(i-1)-j)) \dots  \notag
\end{multline}

\noindent for $1\leq i\leq a$ shows that $H^{1}(S;\OS(-aL)(-j))=0$
iff $j>(a-1)(k-2)$ or $j=0,a=0,1$ , proving iv). The sequences

\begin{multline}
0 \longrightarrow H^{0}(S;\OS(-iC)(-j)) \longrightarrow H^{0}(S;\OS(-(i-1)C)(-j))
\longrightarrow H^{0}(C;\mO_{C}(-j)) \\
\longrightarrow H^{1}(S;\OS(-iC)(-j)) \longrightarrow H^{1}(S;\OS(-(i-1)C)(-j))
\longrightarrow \dots  \notag
\end{multline}

\noindent for $1\leq i \leq b$ imply that $H^{1}(S;\OS(-bC)(-j))=0$
exactly when $j>0,b\geq 0$ and $j=0, b=0,1$ , proving v). The
cohomology sequence of

$$
\sexs{\Op{3}(j-k)}{\Op{3}(j)}{\OS(j)}
$$

\noindent gives $h^{0}(S;\OS(j))=\binom{j+3}{3}-\binom{j-k+3}{3}.$
Similarly, $h^{0}(C;\mO_{C}(j))=\binom{j+2}{2}-\binom{j-k+3}{2}.$
The sequences

\begin{multline}
0 \longrightarrow H^{0}(S;\OS((i-1)C)(j)) \longrightarrow H^{0}(S;\OS(iC)(j))
\longrightarrow H^{0}(C;\mO_{C}(j)) \\
\longrightarrow H^{1}(S;\OS((i-1)C)(j))\longrightarrow \dots  \notag
\end{multline}

\noindent $1\leq i\leq b$ and the vanishing $H^{1}(S;\OS((i-1)C)(j))\cong
H^{1}(S;\OS(-(i-1)C)(-j+k-4))=0 \text{ for } j>k-4$ (by part v))
give

$$
h^{0}(S;\OS(bC)(j))=h^{0}(S;\OS(j))+bh^{0}(C;\mO_{C}(j))
$$

\noindent and vi) follows.
\end{proof}

\bigskip

\begin{Lem} \label{Lem:sections}
If $a \text{ or } b<0,\quad h^{0}(S;\OS(aL+bC))=0.$ When $a,b\geq
0$:

\noindent For $b\geq a\geq k-3$ or $a\geq b$ and $b(k-1)-a(k-2)\geq 0$,

$$
h^{0}(S;\OS(aL+bC))=(k-1)ab-\frac{(k-2)}{2}a^{2}-\frac{(k-4)}{2}(a+(k-1)b)
+\binom{k-1}{3}+1.
$$

\noindent For $a\geq b$ and $j_{0}$ the largest integer between 0 and $a-b$ such that
$b-j_{0}(k-2)\geq 0$,

$$
h^{0}(S;\OS(aL+bC))=\binom{b+3}{3}-\binom{b-k+3}{3}+(b+1)j_{0}-(k-2)\binom{j_{0}+1}{2}.
$$

\noindent For $b\geq a, a\leq k-2$,

$$
h^{0}(S;\OS(aL+bC))=\binom{a+2}{2}[b-\frac{2a}{3}+1].
$$
\end{Lem}

\begin{proof}
For $a,b<0$ , it is clear that $h^{0}(S;\OS(aL+bC))=0$. If $b\geq 0>a$
, $h^{0}(S;\OS(aL+bC))=h^{0}(S;\OS((a-b)L)(b)=0$ by
Lemma~\ref{Lem:Main Lemma}i. The case $a\geq 0>b$ is handled in
the same way.

Assume $b\geq a\geq k-3.$ Then
$h^{0}(S;\OS(aL+bC))=h^{0}(S;\OS((b-a)C)(a))$ and
by Lemma~\ref{Lem:Main Lemma}vi this is

$$
\binom{a+3}{3}-\binom{a-k+3}{3}+(b-a)[\binom{a+2}{2}-\binom{a-k+3}{2}]
$$

\noindent which is easily shown to equal our formula.

If $b\geq a$ and $a\leq k-2$, we use the mapping $S\overset{\pi}{\longrightarrow}\PS{1}$
and $H^{0}(S;\OS((b-a)C)(a))\cong H^{0}(\PS{1};\pi_{*}\OS((b-a)C)(a)).$
To calculate the direct image sheaf, let the homogeneous
coordinates of $\PS{3}\equiv \mP V$ be chosen so that the line L is given by
$x_{2}=0,x_{3}=0.$
Set $W=\{\xi=\xi_{2}x_{2}+\xi_{3}x_{3}\}\subset H^{0}(S;\OS(1))\cong V^{*}.$
If $S \hookrightarrow \PS{3}$ is defined in by $g=0$ then
$g=x_{2}g_{2}+x_{3}g_{3}$, for $g_{2}$, $g_{3}$ of degree $k-1.$
$S\cap H_{\xi}=L+C_{\xi}$ where $C_{\xi}$ is defined by
$g_{\xi}\equiv \xi_{3}g_{2}-\xi_{2}g_{3}=0.$
Then $S\overset{\pi}{\longrightarrow}\mP W^{*}$ is given by
$\pi(p)=\{\xi \in W \mid g_{\xi}(p)=0\}$ and $C_{\xi}$ is the fiber over
$\xi \in \PS{1}.$ Let C be a fixed fiber defined by $t=0$, for t a
coordinate on $\PS{1}.$ Then the isomorphism
$H^{0}(C_{\xi};\mO_{C_{\xi}}((b-a)C)(a))\cong
H^{0}(C_{\xi};\mO_{C_{\xi}}(a))$ is given by $t^{b-a}s\mapsto s$
for $s\in H^{0}(C_{\xi};\mO_{C_{\xi}}(a)).$ For $Ann\xi \equiv \{x\in V \mid \xi
(x)=0\}$, $C_{\xi}$ is a curve in $\mP Ann\xi$ of degree k-1 and
we have a restriction isomorphism $H^{0}(\mP Ann\xi;\mO_{\mP
Ann\xi}(a))\cong H^{0}(C_{\xi};\mO_{C_{\xi}}(a))$ when $a\leq k-2.$
Writing $V^{*}=U\oplus W$ for $U=span\{x_{0},x_{1}\}$,

\begin{align}
H^{0}(\mP Ann\xi;\mO_{\mP Ann\xi}(a))&\cong \Sym^{a}(V^{*}/\cnums \xi)
\notag \\
\quad &\cong \bigoplus_{i=0}^{a} \Sym^{i}U \otimes
\Sym^{a-i}(W/\cnums \xi) \notag
\end{align}

\noindent which gives

\begin{align}
\pi_{*}\OS((b-a)C)(a)&\cong \Op{1}(b-a)\otimes \bigoplus_{i=0}^{a}
\Sym^{i}U\otimes \Op{1}(a-i) \notag \\
\quad &\cong \bigoplus_{i=0}^{a}{\Op{1}(b-i)}^{\oplus i+1}. \notag
\end{align}

\noindent Therefore

\begin{align}
h^{0}(S;\OS((b-a)C)(a)) &=\sum_{i=0}^{a} (i+1)(b-i+1) \notag \\
\quad &=\binom{a+2}{2}[b-\frac{2a}{3}+1]. \notag
\end{align}

Now assume $a\geq b$ and $b(k-1)-a(k-2)\geq 0.$
Then $h^{0}(S;\OS(aL+bC))=h^{0}(S;\OS((a-b)L)(b))$ and consider the cohomology of the
sequences

$$
\sexs{\OS((j-1)L)(b)}{\OS(jL)(b)}{\mO_{L}(b-j(k-2))}
$$

\noindent for $1\leq j\leq a-b.$ Since $b(k-1)-a(k-2)\geq 0$, $b-j(k-2)\geq 0$ for all
j and $h^{1}(S;\OS((j-1)L)(b))=h^{1}(S;\OS((1-j)L)(k-4-b))=0$
by Lemma~\ref{Lem:Main Lemma}iv because
$b-k+4>(k-2)(a-b-2)$ is $b(k-1)-a(k-2)+k>0.$ This gives

$$
h^{0}(S;\OS((a-b)L)(b))=\binom{b+3}{3}-\binom{b-k+3}{3}+(a-b)(b+1)-(k-2)\binom{a-b+1}{2}
$$

\noindent which is equivalent to our formula.

If $a\geq b$ but $b-j(k-2)<0$ for some $1\leq j\leq a-b$ , the above
argument is easily adjusted to give our result.
\end{proof}

\bigskip

\section{Examples of Rank 2 Bundles} \label{S:examples of rank
2}

Let $S\in \PS{3}$ be a smooth surface of degree $k=2\nu +c_{1}$
containing a line L, $\nu \in \znums_{+}, c_{1}=0 \text{ or } -1$,
and $\mL$ a line bundle on S determined by $\mL=\OS(-aL-bC)(\nu+c_{1})$
where  $a,b\in \znums.$ This gives
$\mL^{*}(\nu+c_{1})=\OS(aL+bC).$ As in section 3, we examine the rank 2
extensions

\begin{equation} \label{E:extension 5}
\begin{CD}
0 @>>> \overset{2}{\oplus}\Op{3}(-\nu) @>> \sigma_{1}\oplus \sigma_{2}> E
@>>> {j_{S}}_{*}\OS(-aL-bC)(\nu+c_{1}) @>>> 0
\end{CD}
\end{equation}

\noindent and determine which divisors $aL+bC$ have the property that
the generic extension (\ref{E:extension 5}) is a stable bundle.
Recall that in the dual sequence, in the case that E is locally
free,

\begin{equation} \label{E:dual equation 5 for E}
\begin{CD}
0 @>>> E^{*} @>>(\sigma_{1}^{t},\sigma_{2}^{t})>
\overset{2}{\oplus}\Op{3}(\nu) @>>\tau_{1}\oplus\tau_{2}>
{j_{S}}_{*}\OS(aL+bC)(\nu) @>>> 0 \\
\end{CD}
\end{equation}

\noindent $\tau \in \overset{2}{\oplus}H^{0}(S;\OS(aL+bC))$ is the
extension class of (\ref{E:extension 5}).

\begin{Thm}  \label{Thm:Main Theorem rank 2}
The generic extensions of the form (\ref{E:extension 5}) with $D=aL+bC$
are stable rank 2 bundles in exactly the following cases:(recall
$k=2\nu +c_{1})$
\begin{enumerate}
\item For $k=2 (\nu =1, c_{1}=0)$, $a=0$ and $b\geq 2$ or vice versa.
Here S is a smooth quadric Q.
Using the bidegree notation for line bundles on Q, either
$\mL=\mO_{Q}(1,1-b)$ for $b\geq 2$ and
$\mL^{*}(\nu+c_{1})=\mO_{Q}(0,b)$ or $\mL=\mO_{Q}(1-b,1)$ and
$\mL^{*}(\nu+c_{1})=\mO_{Q}(b,0).$ $c_{2}(E)=b-1.$
\item For $k\geq 3$, $a=0$ and $b>\nu+c_{1}.$ $\mL^{*}(\nu+c_{1})=\OS(bC)$
and $\mL=\OS(-bC)(\nu+c_{1}).$
$c_{2}(E)=b(2\nu+c_{1}-1) -\nu(\nu+c_{1})=b(k-1)-(k^{2}-c_{1}^{2})/4.$
\end{enumerate}
Note that, if $k=2$ in statement (2), statement (1) results. No
other values of $a$ and $b$ produce stable bundles.
\end{Thm}

\begin{proof}
From Section~\ref{S:sb,r2} we know
that k=1 can not occur and that the generic extension E is a stable
bundle iff $\OS(D)$ is globally generated, $D^{2}=0$ , and $h^{0}(S;\mL)=0.$
$D^{2}=0$ gives $a^{2}(2-k)+2ab(k-1)=0$ and so

\begin{equation} \label{E:a and b 1}
a=0 \text{ or } 2b(k-1)-a(k-2)=0.
\end{equation}

\noindent Since $\mO_{L}(D)= \mO_{L}(aL^{2}+bC\cdot
L)=\mO_{L}(b(k-1)-a(k-2))$ is globally generated,

\begin{equation} \label{E:a and b 2}
b(k-1)-a(k-2)\geq 0.
\end{equation}

\noindent Since $\mO_{C}(D)$ is globally generated,
$0\leq \text{ deg}\mO_{C}(D)=C\cdot (aL+bC)$ , i.e.

\begin{equation} \label{E:a and b 3}
a(k-1)\geq 0.
\end{equation}

\noindent If k=2, then $\nu=1$ , $c_{1}=0$ , S=Q, and L and C are
lines from the two pencils of lines on Q. The equations give that
either $a$ or $b=0$ and the other is non-negative. We can assume
$a=0.$ Then $\mL \cong \mO_{Q}(1,1-b)$ and $h^{0}(Q;\mL)=0$ implies $b\geq 2.$
Now $c_{2}(E)=b-1$ follows from (\ref{E:second chern class of E}).

If $k\geq 3$ , $a=0$ because otherwise (\ref{E:a and b 1}) and
(\ref{E:a and b 2}) give $b(k-1)\geq a(k-2)=2b(k-1)$ and so $0\geq b(k-1)$
i.e. $b\leq 0.$ Now (\ref{E:a and b 2}) and (\ref{E:a and b 3})
give $a=0=b.$ But then $\mL=\OS(\nu+c_{1})$ has non-zero global sections. Now
$\mL=\OS(-bC)(\nu+c_{1})$ will have $h^{0}=0$ iff $b>\nu +c_{1}$ by
Lemma~\ref{Lem:Main Lemma} ii.
\end{proof}

\pagebreak

\begin{Prop} \label{Prop:properties of rank 2 bundles}
Let $E\longrightarrow \PS{3}$ be a rank 2 stable bundle of the
type constructed in Theorem~\ref{Thm:Main Theorem rank 2}. Then \\
i) For $l\geq -c_{1}-4, \quad H^{3}(\Op{3};E(l))=0.$ \\
ii) For $l>\nu-4, \quad H^{2}(\Op{3};E(l))=0.$ \\
iii) For $l>b(k-1)-\nu -c_{1}-2, \quad H^{1}(\Op{3};E(l))=0.$ \\
iv) $E(l)$ is globally generated iff $l\geq b(k-1)-\nu -c_{1}.$ \\
v) The line $L\subset S$ is a jumping line of E of jump size $m=b(k-1)-\nu$
i.e. $E_{L}\cong \mO_{L}(m)\oplus \mO_{L}(-m+c_{1}).$
\end{Prop}

\begin{proof}
$H^{3}(\Op{3};E(l))\cong H^{0}(\Op{3};E(-l-c_{1}-4))=0$ for $-l-c_{1}-4\leq 0$
because E is stable so i) holds. The cohomology sequence of

$$
\sexs{\overset{2}{\oplus} \Op{3}(l-\nu)}{E(l)}{{j_{S}}_{*}\OS(-bC)(\nu+c_{1}+l)}
$$

\noindent gives $H^{2}(\Op{3};E(l)) \cong
H^{2}(S;\OS(-bC)(\nu+c_{1}+l))$ when $l>\nu-4.$ $H^{2}(S;\OS(-bC)(\nu+c_{1}+l))\cong
H^{0}(S;\OS(bC)(-\nu-c_{1}-l+k-4))^{*}\cong H^{0}(S;\OS(-bL)(\nu -4-l+b))^{*}=0$
iff $l>\nu-4$ by Lemma~\ref{Lem:Main Lemma}i. This gives ii). The
sequence also gives $H^{1}(\Op{3};E(l)) \cong
H^{1}(S;\OS(-bC)(\nu+c_{1}+l))\cong H^{1}(S;\OS(bC)(\nu-4-l))^{*}
\cong H^{1}(S;\OS(-bL)(b+\nu-4-l))^{*} =0$ for $l+4-b-\nu>(k-2)(b-1)$
, that is, $l>b(k-1)-\nu -c_{1}-2$ by Lemma~\ref{Lem:Main
Lemma}iv. This proves iii). From

$$
\sexs{\overset{2}{\oplus} H^{0}(\PS{3};\Op{3}(l-\nu))}{H^{0}(\PS{3};E(l))}{H^{0}(S;\OS
(-bC)(\nu +c_{1}+l))}
$$

\noindent one sees that E(l) is globally generated iff 1) $l\geq \nu$
and 2) $\OS(-bC)(\nu +c_{1}+l) \cong \OS(bL)(\nu +c_{1}+l-b)$ is
globally generated. A necessary condition for 2) is $i\equiv \nu
+c_{1}+l-b\geq 0.$ (Lemma~\ref{Lem:Main Lemma}ii). The cohomology
sequences of

$$
\sexs{\OS((j-1)L)(i)}{\OS(jL)(i)}{\mO_{L}(jL)(i)\cong
\mO_{L}(i-(k-2)j)}
$$

\noindent for j=1 to b show that $i-(k-2)b\geq 0$ is also
necessary. It is also sufficient because $H^{1}(S;\OS((j-1)L)(i))
\cong H^{1}(S;\OS(-(j-1)L)(-i+k-4))^{*}=0$ j=1 to b for
$i-k+4>(b-2)(k-2)$ , that is, $i>b(k-2)-k$ by Lemma~\ref{Lem:Main
Lemma}iv. Thus $l\geq b(k-1)-\nu-c_{1}$ is necessary and
sufficient for 2). Note that $b(k-1)-\nu-c_{1}
\geq (\nu + c_{1}+1)(2\nu + c_{1}-1)-\nu-c_{1}=2\nu^{2}+3c_{1}\nu -1-2c_{1}\geq \nu$
for $k\geq 2.$ Therefore $l\geq b(k-1)-\nu -c_{1}$ is a necessary
and sufficient condition for E(l) to be globally generated.

To examine L as a jumping line of E express $E_{L}=\mO_{L}(m)\oplus \mO_{L}(-m+c_{1})$
for some $m\geq 0$ and restrict (\ref{E:extension 5}) to L to get

$$
\sexs{im({\sigma_{1}}_{L}\oplus {\sigma_{2}}_{L})}{\mO_{L}(m)\oplus
\mO_{L}(-m+c_{1})}{\mO_{L}(\nu +c_{1}-b(k-1))}.
$$

\noindent Because $\nu+c_{1}-b(k-1)<0$ , it is clear that $-m+c_{1}=\nu+c_{1}-b(k-1)$
which gives the result.
\end{proof}

To make some observations about moduli, let $\mM$ be the moduli
space of S-equivalence classes of semi-stable rank two sheaves on $\PS{3}$
 with fixed chern classes $c_{1}=0 \text{ or } -1, c_{2} \text{, and } c_{3}=0$
 , a projective scheme containing the stable rank two bundles as
 an open subset. For E a rank two stable bundle, the Zariski
 tangent space of $\mM$ at E is

\begin{equation} \label{E:Zariski tangent space}
T{{\mM}}_{E} \cong H^{1}(\PS{3};\sheafend(E))
\end{equation}

\noindent and one knows that

\begin{equation} \label{E:dimension of moduli space}
h^{1}(\PS{3};\sheafend(E)) \geq dim_{E}{\mM} \geq
h^{1}(\PS{3};\sheafend(E))- h^{2}(\PS{3};\sheafend(E))
\end{equation}

\noindent and $h^{2}(\PS{3};{\sheafend}(E))=0$ implies that $\mM$ is
smooth at E \cite[Sect. 4.5]{Huybrechts-Lehn}. From Riemann-Roch,

\begin{equation} \label{E:expected dimension, rank 2}
h^{1}(\PS{3};{\sheafend}(E))-h^{2}(\PS{3};{\sheafend}(E))=8c_{2}(E)+2c_{1}-3.
\end{equation}

We want to count the parameters of our construction.
Note that the basic sequences (\ref{E:extension 5}) and (\ref{E:dual
equation 5 for E}) or, more generally, (\ref{E:sequence for E}) and
(\ref{E:dual equation for E}) are dual to one another.
Also note that, when $k\geq 3$, the isomorphism class of
$\mL=\OS(-bC)(\nu+c_{1})\cong \OS(bL)(\nu+c_{1}-b)$ is determined by the line L
because two lines on S (or integer multiples of lines)
can not be linearly equivalent (or even
homologically equivalent): If $L,L^{\prime}\subset S$ are homologically
equivalent, $L^{2}=L\cdot L^{\prime}\geq 0$; but we know $L^{2}=-(k-2)<0.$
Also when $k\geq 3$, S can contain only a finite number of lines. To
see this, let G be the grassmannian of lines and $\mS$ the universal
sub-bundle over G. Then the degree k polynomial g defining S can
be viewed as a global section of $\Sym^{k}(\mS^{*})$ whose zeroes
are the lines contained in S. The zero set of g is either finite
or of positive dimension. In the latter case, since $Pic(S)$ is
discrete, there are linearly equivalent lines on S, a
contradiction.

For fixed $c_{1}$, $\nu$, $a=0$, and b, $(S,L,\tau)$ defines $(E,\sigma)$
and the function $(S,L,\tau)\rightarrow (E,\sigma)$ is injective
but not a priori surjective, as we explain.
From Theorem~\ref{Thm:Main Theorem rank 2}, E has the form

\begin{equation} \label{E:original}
\begin{CD}
0 @>>> \overset{2}{\oplus}\Op{3}(-\nu) @>> \sigma_{1}\oplus \sigma_{2}> E
@>>> {j_{S}}_{*}\OS(-bC)(\nu+c_{1}) @>>> 0.
\end{CD}
\end{equation}

\noindent Choose a different $\bar{\sigma}\in
\overset{2}{\oplus}H^{0}(\PS{3};E(\nu))$; this produces another
sequence

\begin{equation} \label{E:different sections}
\begin{CD}
0 @>>> \overset{2}{\oplus}\Op{3}(-\nu) @>> \bar{\sigma_{1}}\oplus \bar{\sigma_{2}}> E
@>>> {j_{\bar{S}}}_{*}\bar{\mL} @>>> 0
\end{CD}
\end{equation}

\noindent for $\bar{S}$ another smooth surface of degree k and $\bar{\mL}$
a line bundle on $\bar{S}.$  Does $\bar{S}$
contain a line and, if so, is $\bar{\mL}$ of the form
$\mO_{\bar{S}}(-b\bar{C})(\nu+c_{1})$?
We show, somewhat surprisingly, that the answer to
both questions is affirmative. Note that these considerations are
relevant only when $b\leq k$ because $b>k$ implies that $h^{0}(\PS{3};E(\nu))=2$
and so $\bar{\sigma}$ differs from $\sigma$ by a basis change.

The cohomology sequences of (\ref{E:original}) and
(\ref{E:different sections}) and Lemma~\ref{Lem:Main Lemma} imply
that $h^{0}(\bar{S};\bar{\mL}(b-\nu -c_{1})=h^{0}(S;\OS(bL))=1$
Therefore $\bar{\mL}(b-\nu-c_{1}) \cong \mO_{\bar{S}}(\bar{D})$ for $\bar{D}$
effective. The chern class formulas (\ref{E:second chern class of E}) and
(\ref{E:chern2}) imply

\begin{align} \label{E:chern number equality}
deg\bar{D}&=\omega_{0}\cdot c_{1}(\mO_{\bar{S}}(\bar{D}))   \\
&=\omega_{0}\cdot c_{1}(\OS(bL)) \notag \\
&=H\cdot bL \notag \\
&=b \notag
\end{align}

\noindent and

\begin{align}
{c_{1}(\bar{\mL})}^{2}&={c_{1}(\OS(-bC)(\nu+c_{1}))}^{2} \notag \\
{(\bar{D}-(b-\nu-c_{1})H)}^{2}&={((\nu+c_{1})H-bC)}^{2} \notag \\
{\bar{D}}^{2}&=-(k-2)b^{2}.  \label{E:chern squared equality}
\end{align}

\noindent Express $\bar{D}=\sum_{i}m_{i}Y_{i}$ for $Y_{i}$ irreducible
curves and $m_{i}\in {\znums}^{+}.$ The genus formula gives

\begin{align}
g_{i}&=1+\frac{1}{2}({Y_{i}}^{2}+(k-4)degY_{i})  \notag \\
{Y_{i}}^{2}&\geq -(k-2)degY_{i}.   \notag
\end{align}

\noindent Now (\ref{E:chern squared equality}) implies

\begin{align}
-(k-2)b^{2}&=\sum_{i}{m_{i}}^{2}{Y_{i}}^{2}+2\sum_{i<j}m_{i}m_{j}Y_{i}\cdot Y_{j}
\notag \\
&\geq -(k-2)\sum_{i}m_{i}^{2}degY_{i}.   \label{E:inequality}
\end{align}

\noindent Using $\sum_{i}m_{i}degY_{i}=b$,

\begin{align}
\sum_{i}m_{i}^{2}degY_{i}-b^{2}&=\sum_{i}m_{i}^{2}degY_{i}-b\sum_{i}m_{i}degY_{i}
 \notag \\
&=\sum_{i}m_{i}degY_{i}(m_{i}-b) \notag \\
&\leq 0  \label{E:inequality 2}
\end{align}

\noindent with equality if and only if there is only one term in
the sum, $b=m_{1}$, and $degY_{1}=1.$ But (\ref{E:inequality}) and
(\ref{E:inequality 2}) show that equality must hold and so $\bar{D}=b\bar{L}$
for $\bar{L}$ a line on $\bar{S}.$ This gives
$\bar{\mL}\cong \mO_{\bar{S}}(-b\bar{C})(\nu+c_{1})$ for $\bar{H}=\bar{L}+\bar{C}$
a hyperplane section of $\bar{S}.$ We have proven that, for $k\geq3$, there is a
1-to-1 correspondence

\begin{equation} \label{E:count 1}
(S,L,\tau)\longleftrightarrow (E,\sigma).
\end{equation}

Note that a linear change in $(\tau_{1},\tau_{2})$ produces
an isomorphic E and a corresponding linear change in
$(\sigma_{1},\sigma_{2}).$ Similarly, a linear change in $(\sigma_{1},\sigma_{2}).$
does not change S or $\mL$ and produces a linear change is $(\tau_{1},\tau_{2}).$
Let $G_{2} \equiv  G_{2}(H^{0}(\PS{3};E(\nu))$
and ${G^{\prime}}_{2} \equiv G_{2}H^{0}(S;\mL^{*}(\nu +c_{1}))$ be
grassmannians.
For $[\sigma]\in G_{2}$ and $[\tau] \in {G^{\prime}}_{2}$, our
1-to-1 correspondence can be refined to:

$$
(S,L,[\tau])\longleftrightarrow (E,[\sigma]).
$$

When k=2, S is a smooth quadric Q with two linear equivalence
classes of lines, $\pm.$ In this case the 1-to-1 correspondence is
$(E,\sigma)\leftrightarrow (Q,\pm,\tau).$

Denote by ${\mathcal Y}$ the subset of $\mM$ consisting of
isomorphism classes of
stable bundles of the form (\ref{E:extension 5}). Define $dim\mY$
as the number of independent parameters determining E
(see Proposition~\ref{Prop:dim} below). In general we
expect $dim_{E}{\mathcal Y} < dim_{E}{\mM}.$ In
Section~\ref{S:scheme} we will discuss what conditions imply that $\mY$
has a natural scheme structure and that $\mY \hookrightarrow \mM$
is a regular map.

\begin{Prop} \label{Prop:dim}
Let $E\longrightarrow \PS{3}$ be a rank 2 stable bundle of the
type constructed in Theorem~\ref{Thm:Main Theorem rank 2}. Let $\mY \subset \mM$
be the set of these bundles. Then

$$
dim \mY=
\begin{cases}
\binom{k+3}{3}+2b-k &\text{ when } b>k\geq 3 \\
21 &\text{ when } k=3,b=3 \\
11 &\text{ when } k=3,b=2 \\
2b+7 &\text{ when } k=2,b\geq 3 \\
5 &\text{ when } k=2,b=2 \\
\binom{k+3}{3}+2b-k-2\binom{k-b+3}{3} &\text{ when } b\leq
k,k\geq4.
\end{cases}
$$

\end{Prop}

\begin{proof}
From the 1-to-1 correspondence (\ref{E:count 1}),

\begin{align}
dim\mY &= dim\{S\}+dim\{\tau\}-dim\{\sigma\}  \notag  \\
&=\binom{k+3}{3}-1-max(k-3,0)+2h^{0}(S;\OS(bC))-2[2+h^{0}(S;\OS(-bC)(k))].
\notag
\end{align}

\noindent Here we have used $dim\{(S,L)\}=dim\{S\}$ when $k\geq 3$ since S
contains at most a finite number of lines.

When $k=2$, i.e. S is a smooth quadric Q, we can
calculate directly that $h^{0}(Q;\mO_{Q}(0,b))=b+1$ and $h^{0}(Q;\mO_{Q}(2,2-b))=3$
when $b=2$ and 0 when $b\geq 3$ which gives the result in this case.

When $k\geq 3$, using Lemma~\ref{Lem:Main Lemma}iii,

$$
dim\mY= \binom{k+3}{3}-k+2b-2h^{0}(S;\OS(-bC)(k)).
$$

\noindent For $b>k$, $h^{0}(S;\OS(-bC)(k))=0$ by
Lemma~\ref{Lem:Main Lemma}ii. For $b\leq k$, use
$h^{0}(S;\OS(-bC)(k))=h^{0}(S;\OS(bL)(k-b))$ and the sequences

$$
\sexs{\OS(j-1)L)(k-b)}{\OS(jL)(k-b)}{\mO_{L}(k-b-(k-2)j)}
$$

\noindent for $1\leq j \leq b.$ For $j\geq 2$ or $j=1$ and either
$k\geq 4$ or $b\geq 3$, $k-b-(k-2)j<0$ and
so the cohomology sequence gives
$h^{0}(S;\OS(bL)(k-b))=h^{0}(S;\OS(k-b))=\binom{k-b+3}{3}.$
The remaining case, k=3, b=2, yields $h^{0}(\OS(2L)(1))=5.$ This
gives our formula when $k\geq 3.$
\end{proof}

\medskip

For the bundles of Theorem~\ref{Thm:Main Theorem rank 2}, the
formula (\ref{E:expected dimension, rank 2}) becomes

\begin{equation} \label{E:ed1}
h^{1}(\PS{3};{\sheafend}(E))-h^{2}(\PS{3};{\sheafend}(E))=8b(k-1)-2k^{2}-3.
\end{equation}

\noindent Therefore by choosing b large compared to $k^{2}$ , one
gets $dim_{E}{\mathcal Y}$ much smaller that $dim_{E}{\mM}$ but,
choosing $b=k+1$ , one gets, for large k, $dim_{E}{\mathcal Y}>$
the expected dimension of $\mM$ at E.

We now obtain an upper bound for $h^{1}(\PS{3};{\sheafend}(E))$ by deriving an
upper bound for $h^{2}(\PS{3};{\sheafend}(E)).$ When k=2 or 3, this will give
$dim_{E}\mM$ exactly. To set up the framework for these
calculations, write out
(\ref{E:extension 5}) and (\ref{E:dual equation 5 for E}) in this
case,

\begin{equation} \label{E:specific extension 5}
\sexs{\overset{2}{\oplus}\Op{3}(-\nu)}{E}{{j_{S}}_{*}\OS (-bC)(\nu
+c_{1})}
\end{equation}

\begin{equation} \label{E:specific dual equation 5 for E}
\sexs{E^{*}}{\overset{2}{\oplus}\Op{3}(\nu)}{{j_{S}}_{*}\OS
(bC)(\nu)}.
\end{equation}

\noindent Tensor (\ref{E:specific dual equation 5 for E}) with E
to get

\begin{equation} \label{E:sexs for End(E)}
\sexs{{\sheafend}(E)}{\overset{2}{\oplus}E(\nu)}{{j_{S}}_{*}E_{S}
(bC)(\nu)}.
\end{equation}

\noindent Tensoring (\ref{E:specific extension 5}) with $\OS$ and
calculating $\OS {\otimes}_{\Op{3}} \OS (-bC)(\nu+c_{1})\cong \OS (-bC)(\nu+c_{1})$
and ${\sheaftor}^{\Op{3}}_{1}(\OS,\OS (-bC)(\nu+c_{1}))=\OS (-bC)(-\nu)$
yields

$$
0\longrightarrow \OS (-bC)(-\nu) \longrightarrow \overset{2}{\oplus}\OS(-\nu)
\longrightarrow E_{S} \longrightarrow \OS (-bC)(\nu+c_{1})
\longrightarrow 0
$$

\noindent which can be written as

\begin{equation} \label{E:lexs for E on S}
0\longrightarrow \OS \longrightarrow \overset{2}{\oplus}\OS(bC)
\longrightarrow E_{S}(bC)(\nu ) \longrightarrow \OS(k)
\longrightarrow 0
\end{equation}

\noindent which can be broken up into two short exact sequences

\begin{multline}   \label{E:sexs for E on S}
0\longrightarrow \OS \longrightarrow \overset{2}{\oplus}\OS(bC)
\longrightarrow {\mathcal K} \longrightarrow 0  \\
0 \longrightarrow {\mathcal K} \longrightarrow E_{S}(bC)(\nu ) \longrightarrow \OS(k)
\longrightarrow 0.
\end{multline}

\noindent The cohomology sequence of (\ref{E:sexs for End(E)}) and
Proposition~\ref{Prop:properties of rank 2 bundles}ii yield

\begin{multline} \label{E:cohomology of End(E)}
0 \longrightarrow H^{0}(\PS{3};\sheafend(E)) \longrightarrow
\overset{2}{\oplus}H^{0}(\PS{3};E(\nu )) \longrightarrow
H^{0}(S;E_{S}(bC)(\nu )) \longrightarrow  \\
H^{1}(\PS{3};\sheafend(E)) \longrightarrow \overset{2}{\oplus}H^{1}(\PS{3};E(\nu ))
\longrightarrow H^{1}(S;E_{S}(bC)(\nu )) \longrightarrow H^{2}(\PS{3};\sheafend(E))
 \longrightarrow 0
\end{multline}

\noindent and therefore $h^{2}(\PS{3};\sheafend(E)) \leq h^{1}(S;E_{S}(bC)(\nu
)).$ From the cohomology of the second sequence in (\ref{E:sexs
for E on S}) , $h^{1}(S;E_{S}(bC)(\nu)) \leq h^{1}(S;{\mathcal K})$
and from the first sequence, $h^{1}(S;{\mathcal K})\leq
2h^{1}(S;\OS(bC)) +h^{2}(S;\OS).$ Note that
$h^{2}(S;\OS)=h^{0}(S;\OS(k-4))= \binom{k-1}{3} =0$ for k=2,3 and
$h^{1}(S;\OS(bC))=h^{1}(S;\OS(-bC)(k-4))=0$ for k=2,3 by
Lemma~\ref{Lem:Main Lemma}v. Therefore by(\ref{E:dimension of
moduli space}) and (\ref{E:ed1}),

\begin{Thm} \label{Thm:moduli in cases k=2 and 3}
Let E be a stable rank 2 bundle as constructed in
Theorem~\ref{Thm:Main Theorem rank 2}. For k=2 or 3,
$H^{2}(\PS{3};\sheafend(E))=0$
and so the moduli space $\mM$ containing E is smooth at E and
$T\mM_{E} =H^{1}(\PS{3};\sheafend(E).$ Its dimension is

\begin{equation}
dim_{E}\mM =
\begin{cases}
8b-11 &k=2, b\geq 2 \\
16b-21 &k=3, b\geq 2.  \notag
\end{cases}
\end{equation}
\end{Thm}

\bigskip

For $k\geq 4,b>k-4$ the inequalities above only give an estimate for
 $h^{1}(\PS{3};\sheafend(E))$ and thus $dim_{E}\mM$.
Riemann-Roch calculates $\chi (S;\OS(bC))=1+\binom{k-1}{3}-(k-4)(k-1)b/2$ and,
since $h^{0}(S;\OS(bC))=b+1$ and $h^{2}(S;\OS(bC))=h^{0}(S;\OS(-bC))(k-4))=0$
(using $b>k-4$ and Lemma~\ref{Lem:Main Lemma}ii),
$h^{1}(S;\OS(bC))=(k-4)(k-1)b/2-\binom{k-1}{3}+b.$ Putting all this
together,

\begin{equation}
h^{2}(\PS{3};\sheafend(E))\leq (k^{2}-5k+6)b-\binom{k-1}{3}.
\end{equation}

\noindent Therefore by (\ref{E:dimension of moduli space}) and
(\ref{E:ed1}), for $k\geq 4$, $b\geq k-4$,

\begin{multline} \label{E:estimate for dimension of moduli space at E}
8(k-1)b-2k^{2}-3 \leq dim_{E}\mM \leq 8(k-1)b-2k^{2}-3+(k^{2}-5k+6)b-\binom{k-1}{3}
\\
=(k^{2}+3k-2)b-(k^{3}+5k^{2}+11k+12)/6.
\end{multline}

\noindent This shows that, for fixed k and large b, the
codimension of $\mY$ in $\mM$ is at least of order $(8k-10)b.$

Returning to the $k=2,3$ cases and comparing $dim\mY$ with $dim_{E}{\mM}$
shows that equality
holds only when k=2, b=2,3 and k=3, b=2. When k=2, b=2
then $c_{1}(E)=0, c_{2}(E)=1$ and these are the null-correlation
bundles classified by Barth \cite{Barth-P2} and Wever \cite{Wever}. The moduli space of
these stable bundles is isomorphic to $\PS{5}-{\Bbb G}(1,3)$ where
${\Bbb G}$ is the grassmannian of lines in $\PS{3}$
\cite[page 266]{Hartshorne-sb}. When k=2,
b=3 then $c_{1}(E)=0, c_{2}(E)=2$ and these stable bundles were
classified and studied in detail by Hartshorne \cite{Hartshorne-sb}. The moduli
space of these bundles is smooth, irreducible, and of dimension
13. When k=3, b=2 then $c_{1}(E)=-1, c_{2}(E)=2$ and these bundles
were analyzed by Hartshorne and Sols \cite{Hartshorne-Sols}. The moduli space of these
stable bundles is smooth, irreducible, and rational of dimension
11.

\bigskip

\section{Scheme Structures Related to the Parameter Space}
\label{S:scheme}

Recall that $\mY$ denotes the set of isomorphism classes of stable
rank two
bundles of the form

\begin{equation} \label{E:sc1}
\begin{CD}
0 @>>> \overset{2}{\oplus}\Op{3}(-\nu) @>> \sigma_{1}\oplus \sigma_{2}> E
@>>> {j_{S}}_{*}\OS(bL)(\nu+c_{1}-b) @>>> 0
\end{CD}
\end{equation}

\noindent where the extension class
$\tau=(\tau_{1},\tau_{2})\in \overset{2}{\oplus}H^{0}(S;\OS(bC))$
appears as a homomorphism in the dual sequence

\begin{equation} \label{E:sc2}
\begin{CD}
0 @>>> E^{*} @>>(\sigma_{1}^{t},\sigma_{2}^{t})>
\overset{2}{\oplus}\Op{3}(\nu) @>>\tau_{1}\oplus\tau_{2}>
{j_{S}}_{*}\OS(-bL)(\nu+b) @>>> 0.
\end{CD}
\end{equation}

\noindent Note that we are using the isomorphism $\OS(bC)\cong \OS(-bL)(b).$ Here
$\nu$, b, and $c_{1}$ are fixed and S, L, $\tau$, E,
and $\sigma$ vary. We have shown that the sequences give a 1-to-1 correspondence
(see Section~\ref{S:examples of rank 2})

\begin{equation} \label{E:sc3}
(E,\sigma_{1},\sigma_{2}) \longleftrightarrow
(S,L,\tau_{1},\tau_{2})
\end{equation}

\noindent which can be refined to

\begin{equation}   \label{E:sc4}
(E,[\sigma_{1},\sigma_{2}])\longleftrightarrow
(S,L,[\tau_{1},\tau_{2}]).
\end{equation}.

We would like to show that $\mY$ has a natural scheme structure
and that there is a regular map $\mY \rightarrow \mM$ into the
full moduli space but this seems to be the case only under certain
circumstances. To discuss the situation we use auxiliary parameter
spaces

$$
\mY_{2}\equiv \{(S,L,\tau)\}
$$

\noindent (where S is a smooth surface of degree $k=2\nu+c_{1}$
containing the line L and $\tau$ globally generates $\OS(bC)$) and

$$
\mY_{1}\equiv \{(S,L,[\tau])\}.
$$

\noindent By (\ref{E:sc4}) there is an bijective function from $\mY_{1}$ to the set
$\{(E,[\sigma])\}.$
We have $\mY = \{E\}$ and the obvious projection
functions

$$
\mY_{2} \overset{\eta_{2}}{\rightarrow} \mY_{1}
\overset{\eta_{1}}{\rightarrow} \mY.
$$

\noindent We will show that $\mY_{2}$ and $\mY_{1}$ have natural
scheme structures and regular maps into $\mM$. Then we will point
out some situations in which these results descend to $\mY$.

Let P be the projective space of surfaces of degree k and G the
grassmannian of lines in $\PS{3}$. Define
$Z\equiv \{(S,L,p)\mid L\subset S,p\in S\}$
,$W\equiv \{ (S,L)\mid L\subset S\}$, and $\pi:Z\rightarrow W$ the
projection. Z and W are clearly projective varieties. Let $Z_{0}\subset Z$
and $W_{0}\subset W$ to be the Zariski open subsets defined by
requiring that S is smooth. We define a line bundle $\mF$ on $Z_{0}$
such that, for all $(S,L)\in W_{0}$, $\mF|_{\pi^{-1}(S,L)} \cong \OS(-bL)(b).$
Actually we define $\mF^{\prime}$ such that
$\mF^{\prime}|_{\pi^{-1}(S,L)} \cong \OS(-bL)$ and then set
$\mF \equiv \mF^{\prime}\otimes {\pi_{3}}^{*}\Op{3}(b).$

For $(S,L,p) \in Z_{0}$
such that $p \notin L$, set $\mF^{\prime}_{(S,L,p)}\equiv \mO_{Z_{0},(S,L,p)}.$
If $p\in L$, we proceed as follows. Let S be defined by $g(x)=0$
so that $S=[g]\in P$ and let L be given by the two linear equations $l_{1}=0$
and $l_{2}=0$ so that $L=[l_{1}\wedge l_{2}] \in G$. For each $(S,L) \in
W_{0}$, $g=l_{1}g_{1}+l_{2}g_{2}$ for $g_{1}$ and $g_{2}$ of degree
$k-1$. Since S is smooth, at least one of $g_{1}$ and $g_{2}$ does
not vanish in a Zariski open neighborhood in S of the given point $p\in L.$
Assume $g_{1}$ never vanishes. On this neighborhood,

$$
l_{1}=-\frac{l_{2}g_{2}}{g_{1}}
$$

\noindent and so the pencil of hyperplane sections of S containing
L, $\{H_{t}\}$, which are defined by $t_{1}l_{1}+t_{2}l_{2}=0$,
can be expressed as

$$
\frac{l_{2}}{g_{1}}(-t_{1}g_{2}+t_{2}g_{1})=0.
$$

\noindent The local equations for L and $C_{t}$ on S are therefore $l_{2}=0$
and $-t_{1}g_{2}+t_{2}g_{1}=0$ respectively.

As $(S,L,x)$ varies in an open neighborhood U of a point
$(S_{0},L_{0},p_{0})$ of $Z_{0}$,
we need to demonstrate that $l_{2}(x)$ above can be chosen as a regular function
of $(S,L,x).$ This requires knowing that $g_{1}(x)$ has no zeroes
on U and so it is sufficient to show that $g_{1}$ is a regular
function of $(S,L,x).$ By a coordinate change we can assume that $L_{0}$
is defined by $x_{2}=0$ and $x_{3}=0$ and so, for $(S,L,x)\in U$,

\begin{align}
l_{1}&=x_{2}+l_{1}^{\prime}(x_{0},x_{1})  \notag \\
l_{2}&=x_{3}+l_{2}^{\prime}(x_{0},x_{1}).  \notag
\end{align}

\noindent Expanding $g(x)=g(x_{0},x_{1},l_{1}-l_{1}^{\prime},x_{3})$
gives $g=l_{1}g_{1}+\tilde{g_{2}}$ for $\tilde{g_{2}}$ not
involving $x_{2}$ in fact

\begin{align}
\tilde{g_{2}}&=g(x_{0},x_{1},-l_{1}^{\prime},x_{3})  \notag \\
g_{1}&=\frac{g(x)-\tilde{g_{2}}(x)}{l_{1}}. \notag
\end{align}

\noindent This shows that $g_{1}$ and $\tilde{g_{2}}$
are regular functions on $U\subset Z_{0}.$
Expanding $\tilde{g_{2}}(x_{0},x_{1},x_{3})=
\tilde{g_{2}}(x_{0},x_{1},l_{2}-l_{2}^{\prime})$ gives

\begin{align}
\tilde{g_{2}}(x_{0},x_{1},-l_{2}^{\prime})&=0  \notag \\
\frac{\tilde{g_{2}}(x_{0},x_{1},x_{3})}{l_{2}}&=g_{2}  \notag
\end{align}

\noindent and shows that $g_{2}$ is also a regular function.

If $p\in L$, define $\mF^{\prime}_{(S,L,p)}\equiv \{f=l_{1}^{b}h \mid
h \text{ is regular on } Z_{0} \text{ at } (S,L,p)\}$. It is clear
that our two definitions of $\mF^{\prime}$ patch together and give
a line bundle in the form of a subsheaf of the sheaf of total
quotient rings on $Z_{0}$ \cite[page 144]{Hartshorne}.
It follows from the definition that
, for all $(S,L)\in W_{0}$, $\mF|_{\pi^{-1}(S,L)} \cong \OS(-bL)(b).$

Since $h^{0}\equiv h^{0}(S;\OS(-bL)(b))=h^{0}(S;\OS(bC))$ is
constant in $(S,L)$, $F\equiv \pi_{*}\mF$ is a vector bundle of
rank $h^{0}$ on $W_{0}.$ The parameter space $\mY_{2}$ is the
Zariski open subset of $F^{\oplus 2}$ consisting of
$\tau=(\tau_{1},\tau_{2})\in H^{0}(S;\OS(-bL)(b))^{\oplus 2}$ such
that $\tau_{1}$ and $\tau_{2}$ generate $\OS(-bL)(b)$. This
defines the structure of $\mY_{2}$ as a variety and hence as a
scheme \cite[Chapter 2, Proposition 2.6]{Hartshorne}.

Applying geometric invariant theory to the quotient

$$
\mY_{2}\overset{\eta_{2}}{\rightarrow}\mY_{1}
$$

\noindent by the reductive group $GL(2,\cnums)$ gives an induced scheme
structure to $\mY_{1}.$ More precisely, let $U_{1}$ be an affine
open subset of $W_{0}$ and let $q:F^{\oplus 2}\rightarrow W_{0}$ be
the bundle projection. Then for $U_{1}$ small enough,
$\tilde{U_{1}}\equiv q^{-1}(U_{1})\cong U_{1}\times \cnums^{2h^{0}}$ is also affine
and these sets cover $F^{\oplus 2}.$
Let $\hat{U_{1}}$ be defined as the orbit space of $\tilde{U_{1}}$
under the action of $GL(2;\cnums).$ By \cite[Theorem 6.3.1]{LePotier},
$\hat{U_{1}}$ has the structure of
an affine scheme and these structures for different $\hat{U_{1}}$
patch together to give a scheme structure to the orbit space of the
$GL(2;\cnums)$-action on $F^{\oplus 2}.$ Because the equations on $F^{\oplus 2}$
defining $\mY_{2}$ as a Zariski open subset are clearly
$GL(2;\cnums)$-invariant, they determine $\mY_{1}$ as a Zariski
open subset of the orbit space of the $GL(2;\cnums)$-action on $F^{\oplus 2}.$

There is a natural regular map $\mY_{2}\rightarrow \mM$ defined by
using the universal property of $\mM$ as follows. We will define a
family of stable rank 2 vector bundles on $\PS{3}$ parameterized by
$\mY_{2}$, that is, a coherent sheaf $\mE$ on $\mY_{2}\times \PS{3}$
such that, for every $(S,L,\tau)\in \mY_{2}$, $E\equiv \mE|_{(S,L,\tau)\times \PS{3}}$
is given by (\ref{E:sc1}). Since these restrictions have the same
Hilbert polynomial, $\mE$ is flat over $\mY_{2}.$ This defines a
unique regular map $\mY_{2}\overset{f_{2}}{\rightarrow}\mM$
sending closed points of $\mY_{2}$ to closed points of $\mM.$

We construct $\mE$ by first defining a coherent sheaf $\mF_{2}$ on
$\mY_{2}\times \PS{3}$ and a sheaf mapping
$\pi_{2}^{*}\Op{3}(\nu)^{\oplus 2} \overset{\phi}{\rightarrow} \mF_{2}$
such that, for each $(S,L,\tau) \in \mY_{2}$, the restriction of $\phi$
to $(S,L,\tau)\times \PS{3}$ is given by (\ref{E:sc2}). Then $\mE$
is defined as the dual of the kernel of $\phi$. The construction
of $\mF_{2}$ and $\phi$ is very similar to that of $\mF$ above and
so is left to the reader.

Because $\mY_{2}\overset{f_{2}}{\rightarrow}\mM$ is constant on
the fibers of $\mY_{2}\overset{\eta_{2}}{\rightarrow}\mY_{1}$, it
induces a regular map $\mY_{1}\overset{f_{1}}{\rightarrow}\mM$.

Note that the above development is much simpler when $k=2$.
Denoting the two types of lines on a smooth quadric $Q$ by $\pm$,
the parameter schemes have the form $\mY_{2}=\{(Q,\pm,\tau)\}$ and
$\mY_{1}=\{(Q,\pm,[\tau])\}.$

When $max(a,b)>k$, $[\sigma]$ is unique, $\mY_{1}=\mY$, and so we
have a regular map of schemes $\mY \rightarrow \mM$ injective on
closed points. When $max(a,b)\leq k$, the fibers of $\mY_{1} \rightarrow \mY$
are open subsets of the grassmannian of two dimensional
subspaces of $H^{0}(\PS{3};E(\nu)).$ It is not clear that the scheme
structure of $\mY_{1}$ descends to $\mY.$

The arguments and results of this section apply equally well to
the examples of stable rank 3 bundles studied in
Section~\ref{S:examples of rank three, part II}.

\bigskip

\section{Stable Bundles of Rank 3 on $\PS{3}$}   \label{S:Stable
Bundles of Rank 3}

Let E be a rank 3 normalized bundle ($c_{1}=0,
-1,\text{ or } -2$) on $\PS{3}$. For $\nu$ large enough, $E(\nu)$ is globally
generated and, using Kleiman transversality, the generic
$\sigma=(\sigma_{1},\sigma_{2},\sigma_{3})
\in \overset{3}{\oplus}H^{0}(\PS{3};E(\nu))$ produces

\begin{equation}     \label{E:sexs for E rank 3}
\begin{CD}
o @>>> \overset{3}{\oplus}\Op{3}(-\nu) @>> \sigma > E
@>>> {j_{S}}_{*}\mL @>>> 0
\end{CD}
\end{equation}

\noindent for
\begin{enumerate}
\item $S\equiv Z_{\sigma_{1}\wedge
\sigma_{2}\wedge \sigma_{3}}\subset
\PS{3}$ a smooth hypersurface of degree $k=3\nu +c_{1}$ and $\mL$
a line bundle on S.
\item $Z_{\sigma_{i}}$ j=1,2,3 zero cycles consisting of
$c_{3}(E(\nu))=c_{3}(E)+\nu c_{2}(E)+\nu^{2}c_{1}+\nu^{3}$ smooth
points.
\item $Z_{\sigma_{i}\wedge \sigma_{j}} \quad i<j$ smooth curves of
degree $c_{2}(E(\nu))=c_{2}(E)+2\nu c_{1}+3\nu^{2}.$
\end{enumerate}

\begin{Prop} \label{Prop:properties of degeneracy curves}
For $V\equiv span\{\sigma_{1},\sigma_{2},\sigma_{3}\}
\subset H^{0}(\PS{3};E(\nu))$ generic as
above, the two-dimensional linear system of curves $Y=Z_{s_{1}\wedge s_{2}}$ for
$s_{1}\wedge s_{2}\in \wedge^{2}V$ on the surface S satisfies \\
i) \quad The curves Y are connected and the generic Y is
smooth.\\
ii) \quad $Y^{2}=c_{3}(E(\nu))$ \quad (intersection on S) \\
iii) \quad $genus(Y)=1+1/2\{c_{3}(E(\nu))+c_{2}(E(\nu))(k-4)\}$
\end{Prop}

\begin{proof}
Set $P\equiv {\Bbb P}V\equiv \PS{2}$ and
$Z\equiv \{(x,s) \in \PS{3}\times P \mid s(x)=0\}.$
For $(x,s)\in Z, x\in S$ and the smoothness of S implies that x
determines $s$ (the subspace of V that vanishes at x is
one-dimensional). It follows that $Z\cong S.$ In

\begin{equation} \label{E:incidence diagram}
\begin{CD}
Z@> \pi_{1} >> S \\
@V \pi_{2} VV  \\
P
\end{CD}
\end{equation}

\noindent the curves Y are $\pi^{-1}_{2}(L)$ for the lines $L\subset P.$
Y connected follows from the Fulton-Hansen connectedness
theorem \cite{Fulton-Hansen}. Note that
$Z_{\sigma_{1}\wedge \sigma _{2}}\cap Z_{\sigma_{2}\wedge \sigma _{3}}=Z_{\sigma_{2}}
$ ($\supset$ is obvious, $\subset$ results from the fact that S is
smooth). An easy local coordinate argument shows that
$Z_{\sigma_{1}\wedge \sigma _{2}}$ and $Z_{\sigma_{2}\wedge \sigma _{3}}$
meet transversely on S at each point of $Z_{\sigma_{2}}.$
Therefore $Y^{2}=Z_{\sigma_{1}\wedge \sigma _{2}}\cdot Z_{\sigma_{2}\wedge \sigma _{3}}
=c_{3}(E(\nu)).$ The genus formula for Y now follows from the usual
genus formula for the curve Y on the surface S and $K_{S}=(k-4)\omega_{0}.$
\end{proof}

\noindent This proposition is a special case of a more general
result \cite{Vitter-dl}.

Applying $Hom_{\Op{3}}(\quad ,\Op{3})$ to (\ref{E:sexs for E rank
3}) gives

\begin{equation} \label{E:dual sequence for E rank 3}
\begin{CD}
0 @>>> E^{*} @>>\sigma^{t} >
\overset{3}{\oplus}\Op{3}(\nu) @>>\tau >
{j_{S}}_{*}{\mathcal L}^{*}(k) @>>> 0 \\
\end{CD}
\end{equation}

\noindent for
$\tau =(\tau_{1},\tau_{2},\tau_{3})\in \overset{3}{\oplus}H^{0}(S;{\mL}^{*}(2\nu+c_{1})).$
Arguing as in the rank 2 case

\begin{equation}
\tau_{i}=(\sigma_{j}\wedge \sigma_{k})_{S}\wedge \quad \text{for(ijk)
an even permutation of (123)}.
\end{equation}

Applying Grothendieck Riemann-Roch to ${j_{S}}_{*}{\mL}$ yields

\begin{equation} \label{E:c2 of E rank 3}
c_{2}(E)=3\nu^{2}+3\nu c_{1}+c_{1}^{2}-c_{1}(\mL)\cdot \omega_{0}
\end{equation}

\begin{equation} \label{E:c3 for E}
c_{3}(E)=(2\nu+c_{1})^{3}-(3\nu +2c_{1})c_{1}(\mL)\cdot
\omega_{0}+{c_{1}(\mL)}^{2}  \quad \text{(intersection on S)}.
\end{equation}

Now assume E is stable. If $\mL$ has the form $\OS(l)$, $h^{0}(\PS{3};E)=0$
implies $l<0.$ Then (\ref{E:c2 of E rank 3}) implies

\bigskip

\begin{it}
\noindent For a fixed stable bundle E and $\nu$ large enough, any
representation of E of the form(\ref{E:sexs for E rank 3}) implies
that $\mL \ncong \OS(l)$ for any l and therefore S belongs to the
Noether-Lefschetz locus.
\end{it}

\bigskip

As in the rank 2 case we want to construct some specific rank 3
stable bundles so reverse the above development, begin with a
given $\nu \in \znums_{+}$, a smooth surface $S\subset \PS{3}$ of
degree $k=3\nu+c_{1}$, a line bundle $\mL$ on S and consider
extensions

\begin{equation} \label{E:extension for E rank 3}
\sexs{\overset{3}{\oplus} \Op{3}(-\nu)}{E}{{j_{S}}_{*}\mL}.
\end{equation}

\noindent They are classified by
$\tau \in \overset{3}{\oplus}H^{0}(S;{\mL}^{*}(2\nu+c_{1})).$
Applying $Hom_{\Op{3}}(\quad ,\Op{3})$ to (\ref{E:extension
for E rank 3}),

\begin{equation} \label{E:dual extension rank 3}
\begin{CD}
0 @>>> E^{*} @>>> \overset{3}{\oplus}\Op{3}(\nu)
@>>\tau > \mL^{*}(k) @>>>
{{\sheafext}^{1}}_{\Op{3}}(E,\Op{3}) @>>> 0.
\end{CD}
\end{equation}

\noindent So E is locally free iff ${\mL}^{*}(2\nu +c_{1})$ is
globally generated by $\tau.$ Therefore

\bigskip

\begin{it}
\noindent The generic extension (\ref{E:extension for E rank 3}) is
locally free iff ${\mL}^{*}(2\nu+c_{1})$ is globally generated
(necessarily by three sections).
\end{it}

\bigskip

A rank 3 reflexive sheaf E on $\PS{3}$ is stable iff $h^{0}(\PS{3};E)=0$
and $h^{0}(\PS{3};E^{*})=0 \quad (c_{1} =0)$, $h^{0}(\PS{3};E^{*}(-1))=0
\quad (c_{1} =-1, -2).$ When $c_{1}=0$ , E is semistable iff $h^{0}(\PS{3};E(-1))=0$
and $h^{0}(\PS{3};E^{*}(-1))=0$ \cite[page 167]{OSS}. Therefore (\ref{E:extension for
E rank 3})
 and (\ref{E:dual extension rank 3}) imply

\bigskip

\begin{it}
\noindent The bundle $E$ of the form (\ref{E:extension for
E rank 3}) is stable iff the following two conditions hold: \\
$A)$ \quad $\nu \geq 1$ and  $h^{0}(S;\mL)=0$ \\
$B)$ \quad $\overset {3}{\oplus}H^{0}(\PS{3};\Op{3}(\nu ))
\overset{\tau}{\longrightarrow} H^{0}(S;{\mL}^{*}(k))$ is
injective $(c_{1}=0)$\\
\quad  $\overset {3}{\oplus}H^{0}(\PS{3};\Op{3}(\nu-1))
\overset{\tau}{\longrightarrow} H^{0}(S;{\mL}^{*}(k-1))$ is
injective $(c_{1}=-1,-2).$ \\
$E$ is semistable $(c_{1}=0)$ iff the following two conditions
hold: \\
$A^{\prime})$ \quad $\nu \geq 0$ and  $h^{0}(S;\mL(-1))=0$ \\
$B^{\prime})$ \quad $\overset {3}{\oplus}H^{0}(\PS{3};\Op{3}(\nu-1 ))
\overset{\tau}{\longrightarrow} H^{0}(S;{\mL}^{*}(k-1))$ is
injective.
\end{it}

\bigskip

\section{Examples of Rank 3 Bundles I} \label{S:examples of rank 3}

First we take $\mL=\OS(-l)$ and examine rank 3 bundles E of the
form

\begin{equation} \label{E:extension using hyperplane bundle}
\begin{CD}
0 @>>> \overset{3}{\oplus}\Op{3}(-\nu ) @>> \sigma > E @>>>
{j_{S}}_{*}\OS(-l) @>>> 0
\end{CD}
\end{equation}

\noindent for $\nu, l \in \znums_{+}$ and show that they are stable for
generic $\sigma.$ The dual sequence is

\begin{equation} \label{E:dual extension using hyperplane bundle}
\begin{CD}
0 @>>> E^{*} @>> \sigma^{t} > \overset{3}{\oplus}\Op{3}(\nu ) @>>
\tau > {j_{S}}_{*}\OS(l+k) @>>> 0.
\end{CD}
\end{equation}

\noindent For E to be locally free, $\tau_{1}, \tau_{2}, \tau_{3}$
must globally generate $\OS(l+2\nu +c_{1}).$ We can identify the $\tau_{i}$
with homogeneous polynomials of degree $l+2\nu+c_{1}$ with no
simultaneous zeroes on S. If g is the degree k homogeneous polynomial
that defines S, (\ref{E:dual extension using hyperplane bundle}) can
be expressed as

\begin{equation} \label{E:2nd seq for E dual}
\begin{CD}
0 @>>> E^{*} @>>> {\Op{3}(\nu )}^{\oplus 3}\oplus \Op{3}(l) @>> \tau
\oplus g > \Op{3}(l+k) @>>> 0
\end{CD}
\end{equation}

\noindent and the condition is that $\tau_{1}, \tau_{2}, \tau_{3}, g$ have
no common zeroes on $\PS{3}$, which holds for generic $\tau$.
To verify stability
condition B) from Section~\ref{S:Stable Bundles of Rank 3}, in the
$c_{1}=0$ case, let $\Gamma_{i} \equiv$ the homogeneous
polynomials of degree i and examine the kernel of

\begin{equation}
\begin{CD}
\Gamma^{\oplus 3}_{\nu } \oplus \Gamma_{l} @>> \tau \oplus g >
\Gamma_{l+k}.   \notag
\end{CD}
\end{equation}

\noindent By \cite[Lemma 3.1]{Dwork}, a pre-Koszul complex graded adaptation of the
Koszul complex, if $\psi =(\psi_{1}, \psi_{2}, \psi_{3}, \psi_{4})$
is in the kernel, $\psi_{i}=\sum_{j}B_{ij}\tau_{j}$ for $\tau _{4}=g$
and $B=(B_{ij})$ a skew-symmetric matrix of homogeneous
polynomials with

$$
deg B_{ij}=
\begin{cases}
l+k-2(l+2\nu +c_{1}) &for 1\leq i,j\leq 3 \\
l+k-(l+2\nu +c_{1})-k &for 1\leq i \leq 3, j=4 \text{ or
vice versa.}
\end{cases}
$$

\noindent In both cases these degrees are negative meaning $B=0$
and $\psi =0.$ E is therefore stable. The stability condition for
the cases $c_{1}=-1 \text{ or } -2$ is checked in the same way.
Taking the dual of (\ref{E:2nd seq for E dual}),

\begin{equation} \label{E:2nd seq for E}
\begin{CD}
0 @>>> \Op{3}(-l-k) @>> (\tau ,g) > {\Op{3}(-\nu)}^{\oplus
3}\oplus \Op{3}(-l) @>> \sigma \oplus h > E @>>>0
\end{CD}
\end{equation}

\noindent where
$\sigma \in \overset{3}{\oplus}H^{0}(\PS{3};E(\nu))$
and $h \in H^{0}(\PS{3};E(l)).$ Note that we may drop the condition
that S be smooth, require only that $\tau$ and g have no common zeroes,
and define the stable bundle directly by (\ref{E:2nd seq for E}). For $\nu$,
 l, and $c_{1}$ fixed, this gives a
one-to-one correspondence

\begin{equation} \label{E:counting}
(\tau , g) \longleftrightarrow (E,\sigma ,h).
\end{equation}

Define $\mY$ to be the set of
stable rank 3 bundles of the form (\ref{E:2nd seq for E}). We will
show that $\mY$ has a natural scheme structure, that the
inclusion $\mY \hookrightarrow \mM$ is a regular map, and that $dim\mY=dim\mM.$
 To begin, define

$$
\mY_{0}\equiv \{(\tau,g)\in \Gamma_{l+k-\nu}^{\oplus 3} \times \Gamma_{k}\mid
\tau_{1}, \tau_{2}, \tau_{3}, g \text{ have
no common zeroes on } \PS{3} \}.
$$

\noindent $\mY_{0}$ is a Zariski open subset of an affine space
and there is a coherent sheaf $\mE$ on $\mY_{0}\times \PS{3}$, flat over
$\mY_{0}$,
such that, for each $(\tau,g)\in \mY_{0}$, $E=\mE|_{(\tau, g)\times \PS{3}}$
is given by (\ref{E:2nd seq for E}). By the universal property of
$\mM$, there is a unique regular map $f_{0}:\mY_{0} \rightarrow \mM$
sending the closed points of $\mY_{0}$ to closed points of $\mM.$

Let $(\tau,g)\in \mY_{0}$ determine E and the sequence (\ref{E:2nd seq for
E})and $(\bar{\tau},\bar{g})$ the sequence for $\bar{E}$. An easy
argument, using $H^{1}(\PS{3};\Op{3}(j))=0$ for all j and the fact
that elements of $\mY_{0}$ have no common zeroes, shows that $E\cong \bar{E}$
if and only if the isomorphism extends to an isomorphism of
sequences

\begin{equation*}
\begin{CD}
0 @>>> \Op{3}(-l-k) @>> (\tau ,g) > {\Op{3}(-\nu)}^{\oplus
3}\oplus \Op{3}(-l) @>> \sigma \oplus h > E @>>>0 \\
@VVV  @VidVV @V{\psi}VV @VVV @VVV \\
0 @>>> \Op{3}(-l-k) @>> (\bar{\tau} ,\bar{g}) > {\Op{3}(-\nu)}^{\oplus
3}\oplus \Op{3}(-l) @>> \bar{\sigma} \oplus \bar{h} > \bar{E} @>>>0
\end{CD}
\end{equation*}

\noindent where

$$
\psi=\left( \begin{matrix}A&v(x) \\
w^{t}(x)&b
\end{matrix} \right)
$$

\noindent for $A\in GL(3;\cnums), b\in {\cnums}^{*}, v(x)\in
\Gamma_{l-\nu}^{\oplus 3},
\text{ and } w(x)\in \Gamma_{\nu-l}^{\oplus 3}.$ The homomorphisms $\psi$
form a Lie group H whose dimension equals
$10+3\binom{|l-\nu|+3}{3}$ if $l\ne \nu$ and
16 if $l=\nu.$ We have shown that $\mY=\mY_{0}/H.$ Using the fact that
the isomorphisms of E are scalar multiples of the identity and
that this multiple is fixed by requiring that the left vertical
homomorphism above is the identity, we can compute

\begin{align} \label{E:dim of Y rank 3}
dim\mY &=dim\{\tau\}+dim\{g\}-dimH     \\
       &=3\binom{l+k-\nu+3}{3} + \binom{k+3}{3}-\begin{cases}
       10+3\binom{|l-\nu|+3}{3} &\text{ for } l\ne \nu  \\
       16 &\text{ for } l=\nu.   \notag
       \end{cases}
\end{align}

To identify the scheme structure of $\mY$ as a Zariski open subset
of a projective scheme $\bar{\mY}$, consider the three cases
 $l>\nu$, $l=\nu$, and $l<\nu.$
When $l>\nu$, the action of H is given by $\bar{\tau}=A\tau+gv$ and
$\bar{g}=bg$. The space of orbits $\bar{\mY}$ is therefore the grassmann
bundle $G_{3}(\mW)$ where $\mW$ is the vector bundle on $\mP \Gamma_{k}$
defined by

$$
\sexs{\mO_{\mP \Gamma_{k}}(-1)\otimes \Gamma_{l-\nu}}{\mO_{\mP \Gamma_{k}}\otimes
\Gamma_{k+l-\nu}}{\mW}
$$

\noindent For $l=\nu$, $H=GL(4;\cnums)$ and $\bar{\mY}=G_{4}(\Gamma_{l}).$ Finally for
$l<\nu$, $\bar{\tau}=A\tau$ and $\bar{g}=w^{t}\tau+bg$ and this
gives $\bar{\mY}=\mmP roj(\Sym \mF)$ for $\mF$ the coherent sheaf over $G_{3}(\Gamma_{k+l-\nu})$
defined by

\begin{equation*}
\begin{CD}
\mK @>\phi >> \mO_{G_{3}(\Gamma_{k+l-\nu})}\otimes \Gamma_{k} @>>> \mF @>>> 0
\end{CD}
\end{equation*}

\noindent where $\mK$ is the vector bundle with fiber ${\Gamma_{\nu -l}}^{\oplus 3}$
associated with the principal frame bundle of the tautological sub-bundle on
$G_{3}(\Gamma_{k+l-\nu})$ by the group action $B(\tau,w)=
(B\tau,{B^{t}}^{-1}w)$ for $B\in GL(3;\cnums)$ and where $\phi(\tau,w)=w^{t}\tau.$
$\bar{\mY}$ is a projective scheme \cite[Chapter II,
Proposition 7.10]{Hartshorne}.

It is clear that $f_{0}$ induces a regular map $f:\mY \rightarrow \mM.$

We now compute $h^{1}(\PS{3};\sheafend(E)).$
The cohomology sequence of (\ref{E:2nd seq for E})
implies $H^{1}(\PS{3};E(j))=0 \quad \forall j$ and
$H^{2}(\PS{3};E(k+l))=0.$ Tensoring (\ref{E:2nd seq for E dual})
with E gives

$$
\sexs{\sheafend(E)}{{E(\nu)}^{\oplus 3}\oplus E(l)}{E(l+k)}
$$

\noindent and therefore

\begin{multline} \label{E:coho seq for End(E)}
0 \longrightarrow H^{0}(\PS{3};\sheafend(E)) \longrightarrow
{H^{0}(\PS{3};E(\nu))}^{\oplus 3} \oplus H^{0}(\PS{3};E(l))
\longrightarrow \\
H^{0}(\PS{3};E(k+l)) \longrightarrow
H^{1}(\PS{3};\sheafend(E)) \longrightarrow 0
\end{multline}

\begin{equation}
H^{2}(\PS{3};\sheafend(E))\cong {H^{2}(\PS{3};E(\nu))}^{\oplus 3}
\oplus H^{2}(\PS{3};E(l)).
\end{equation}

\noindent It follows from (\ref{E:extension using hyperplane
bundle}) or (\ref{E:2nd seq for E}) that

$$
h^{0}(\PS{3};E(\nu))=
\begin{cases}
3 &l>\nu \\
3+\binom{\nu -l+3}{3} &l\leq \nu
\end{cases}
$$

\medskip

$$
h^{0}(\PS{3};E(l))=
\begin{cases}
1+3\binom{l-\nu+3}{3} &l\geq \nu  \\
1 &l<\nu
\end{cases}
$$

\medskip

$$
h^{0}(\PS{3};E(k+l))= 3\binom{l+k-\nu+3}{3} + \binom{k+3}{3}
-1
$$

\noindent and so (\ref{E:dim of Y rank 3}) and (\ref{E:coho seq
for End(E)}) imply

\begin{Thm}
Let E be a stable rank 3 bundle on $\PS{3}$ of the
form

\begin{equation}
\begin{CD}
0 @>>> \Op{3}(-l-k) @>> (\tau ,g) > {\Op{3}(-\nu)}^{\oplus
3}\oplus \Op{3}(-l) @>> \sigma \oplus h > E @>>>0. \notag
\end{CD}
\end{equation}

\noindent If $\mY \overset{f}{\hookrightarrow} \mM$ is the set of these bundles, then

\begin{align} \label{E:h1endE}
h^{1}(\PS{3};\sheafend(E))&=dim\mY =dim_{E}\mM    \\
\quad &=3\binom{l+k-\nu+3}{3} + \binom{k+3}{3}-\begin{cases}
       10+3\binom{|l-\nu|+3}{3} &\text{ for } l\ne \nu  \notag \\
       16 &\text{ for } l=\nu.   \notag
       \end{cases}
\end{align}

\noindent $\mY$ is an open subscheme of $\mM$ and $\mM$ is
smooth at E.
\end{Thm}

\bigskip

The second chern class of these bundles is (from (\ref{E:c2 of E
rank 3}))

$$
c_{2}(E)=3\nu^{2}+3\nu c_{1}
+c_{1}^{2}+lk=\frac{1}{3}[k^{2}+c_{1}k+c_{1}^{2}]+lk.
$$

\noindent From Riemann-Roch one gets

\begin{equation}   \label{E:exp dim of moduli rank 3}
h^{1}(\PS{3};\sheafend(E))-h^{2}(\PS{3};\sheafend(E))=12c_{2}(E)-4c_{1}^{2}-8.
\end{equation}

\noindent Combining the above two equations with (\ref{E:h1endE}) gives

\begin{multline}
h^{2}(\PS{3};\sheafend(E))=3\binom{l+k-\nu+3}{3} +
\binom{k+3}{3}-4k[k+c_{1}+3l] \notag \\
-\begin{cases}2+3\binom{|l-\nu|+3}{3} &\text{ if } l\ne \nu \\
8 &\text{ if } l=\nu.
\end{cases}
\end{multline}

Returning to the scheme structure of $\mY$ and the regular map
$\mY \overset{f}{\hookrightarrow} \mM$,
it is tempting to suppose that the closure of $\mY$ in $\mM$
is the projective scheme $\bar{\mY}$ given above. This is not the
case unless $k=l=\nu=1$ in which case $\mY$ is a point and $E\cong T\PS{3}(-2).$
When $\nu>1$, there are points $[\tau, g]\in \bar{\mY} \smallsetminus \mY$ which
correspond to reflexive sheaves of rank three which are unstable.
More precisely, let $Z\subset \PS{3}$ be the subscheme defined by
the vanishing of $(\tau, g)$ and define E by (\ref{E:2nd seq for
E}). For $codimZ\geq 2$, E is torsion free, locally free on $\PS{3}\smallsetminus Z.$
Apply $\sheafhom_{\Op{3}}(\quad,\Op{3})$ to get

$$
\begin{CD}
0 @>>> E^{*} @>>> {\Op{3}(\nu )}^{\oplus 3}\oplus \Op{3}(l) @>> \tau
\oplus g > \mI(l+k) @>>> 0
\end{CD}
$$

$$
\sheafext_{\Op{3}}^{1}(E,\Op{3})\cong \mO_{Z}(k+l)
$$

\noindent for $\mI$ the ideal sheaf of Z. $E^{*}$ is reflexive
(the dual of any coherent sheaf is reflexive) and the sequence
exhibits $E^{*}$ as a second syzygy sheaf. If Z is a 0-dimensional
locally complete intersection, taking the dual again shows that E
is reflexive.
Now take $\tau \equiv (x_{0}^{k+l-\nu-1}x_{3},x_{1}^{k+l-\nu-1}x_{2},
x_{0}^{k+l-\nu-1}x_{2})$ and $g\equiv \sum_{i=0}^{3} x_{i}^{k}.$
Then Z is a zero-dimensional locally complete intersection in $\PS{3}.$
Define $\tilde{f}\equiv (x_{2},0,-x_{3},0)$,
a section of ${\Op{3}(1)}^{\oplus 3}\oplus \Op{3}(l-\nu +1).$
Because $\tilde{f}$ is in the kernel of $\tau \oplus g$, it
defines a section f of $E^{*}(-(\nu -1)).$ This implies that E is unstable
when $\nu>1.$

The family (\ref{E:2nd seq for E}) also contains E which fail to
be torsion free: if $l-\nu\geq 1$ and
$\tau = (x_{1}^{l-\nu}g, x_{2}^{l-\nu}g, x_{3}^{l-\nu}g)$ then E
has torsion. When $l=\nu$ and $k\geq 2$, it is also easy to
construct E with torsion.

\bigskip

\section{Examples of Rank Three Bundles II} \label{S:examples of rank three, part II}

We now construct examples from surfaces $S\hookrightarrow \PS{3}$
containing a line L, choosing a line bundle on S of the
form $\mL=\OS(-aL-bC)(2\nu +c_{1})$ for $a,b\in \znums.$
Then $\mL^{*}(2\nu +c_{1})=\OS(aL+bC).$ We analyze rank 3 extensions

\begin{equation} \label{E:E from L rank 3}
\sexs{\overset{3}{\oplus}\Op{3}(-\nu )}{E}{{j_{S}}_{*}\OS(-aL-bC)(2\nu
+c_{1})}
\end{equation}

\noindent and determine the divisors $aL+bC$ for which the generic
extension is a stable bundle. The dual sequence is

\begin{equation} \label{E:dual E from L rank 3}
\begin{CD}
0 @>>> E^{*} @>>> \overset{3}{\oplus}\Op{3}(\nu ) @> \tau >>
{j_{S}}_{*}\OS(aL+bC)(\nu) @>>> 0
\end{CD}
\end{equation}

\noindent and $\tau \in \overset{3}{\oplus}H^{0}(S;\OS(aL+bC))$ is the
extension class.

\bigskip

\begin{Thm} \label{Thm:main theorem rank 3}
The generic extension of the form (\ref{E:E from L rank 3})
is a stable rank 3 bundle in the following cases:(recall
$k=3\nu +c_{1}$)
\begin{enumerate}
\item $k=1, \nu =1, c_{1}=-2, S= \text{ hyperplane } H.$ There
is no curve C, $a>0$, and $\mL =\mO_{H}(-a).$ E is a special case
of (\ref{E:2nd seq for E}). $c_{2}(E)=a+1$ and $c_{3}(E)=a^{2}-a.$
\item $k=2, \nu =1, c_{1}=-1, S=\text{ smooth quadric } Q$ , and L
and C belong to the two pencils of lines on Q. Using the bidegree
notation for line bundles on Q,
$\mL^{*}(2\nu +c_{1})=\mO_{Q}(a,b), \mL=\mO_{Q}(1-a,1-b)$ for
$a,b\geq 0, max(a,b)\geq 2.$ $c_{2}(E)=a+b-1$ and $c_{3}(E)=2ab-a-b+1.$
\item $k\geq 3$ and $a>b\geq \frac{(k-2)}{(k-1)}a>0$ (which
implies $a\geq k-1, b\geq k-2$) ---except for the case k=3, $a=2$,
$b=1$. $c_{2}(E)=a+b(k-1)-[k^{2}-c_{1}^{2}]/3$ and
$c_{3}(E)=2ab(k-1)-a^{2}(k-2)-(a+(k-1)b)(k-c_{1})/3
+{(k-c_{1})}^{2}(2k+c_{1})/27.$
\item $k\geq 3, b\geq a>\nu /2 (c_{1}=0), (\nu -1)/2 (c_{1}=-1,
-2)$, and $b>2\nu +c_{1}.$ The chern classes of E are as in case
3.
\end{enumerate}
\noindent When $k\geq 3, b\geq a, b>2\nu +c_{1}$ but
$a\leq \nu/2 (c_{1}=0), (\nu -1)/2 (c_{1}=-1,-2)$ the generic
extension is locally free but it is not known if it is stable. No
other values of $a$ and $b$ produce stable bundles.

\end{Thm}

\begin{proof}
The conditions from Section~\ref{S:Stable Bundles of Rank 3} for E
to be locally free and stable are applied to $\mL=\OS(-D)(2\nu +c_{1})$
and $\mL^{*}(2\nu +c_{1})=\OS(D)$ for $D=aL+bC.$ Case (1) follows from the
analysis of the bundles (\ref{E:extension using hyperplane
bundle}). If $\OS(aL+bC)$ is globally generated, $\mO_{L}(aL+bC)=\mO_{L}(b(k-1)-a(k-2))$
is globally generated and so

\begin{equation} \label{E:eq1}
b(k-1)-a(k-2)\geq 0
\end{equation}

\noindent and $\mO_{C}(aL+bC)$ is globally generated and of degree
$C\cdot (aL+bC)=(k-1)a$ so

\begin{equation} \label{E:eq2}
a(k-1)\geq 0.
\end{equation}

\noindent The condition that $h^{0}(\mL)=0$ is, expressing

\begin{align}
\mL &\cong \OS(-(a-b)L)(2\nu+c_{1}-b) \quad \text{for } a\geq b  \notag \\
\quad &\cong \OS(-(b-a)C)(2\nu+c_{1}-a) \quad \text{for } b\geq a
\notag
\end{align}

\noindent and applying Lemma~\ref{Lem:Main Lemma}, equivalent to

\begin{equation} \label{E:eq3}
max(a,b)> 2\nu +c_{1}.
\end{equation}

For $k\geq 3$ and $a>b$ (case 3),
$1\leq a-b\leq a-\frac{k-2}{k-1}a=\frac{a}{k-1}$ and so $a\geq k-1$
and $b\geq k-2.$ Condition (\ref{E:eq3}) is satisfied except for
the case k=3, a=2, b=1. To show that $\OS(aL+bC)$ is globally
generated first note that $\OS(aL+bC)=\OS((a-b)L)(b)$ is clearly
globally generated on $S\smallsetminus L.$ For $1\leq j\leq a-b$, consider the
sequences

\begin{equation} \label{E:eq4}
\sexs{\OS((j-1)L)(b)}{\OS(jL)(b)}{\mO_{L}(b-(k-2)j)}
\end{equation}

\noindent and note that $b-(k-2)j\geq b-(k-2)(a-b)=-(k-2)a+(k-1)b\geq 0$
so that $\mO_{L}(b-(k-2)j)$ is globally generated. The cohomology
sequences now show that $\OS(aL+bC)$ is globally generated because
$H^{1}(S;\OS((j-1)L)(b))\cong
{H^{1}(S;\OS(-(j-1)L)(-b+k-4))}^{*}=0$ by Lemma~\ref{Lem:Main
Lemma}iv since $b-k+4>(k-2)(a-b-2) i.e. (k-1)b-(k-2)a +k>0.$ Now
it will be verified that
$\overset{3}{\oplus}H^{0}(\PS{3};\Op{3}(\nu )) \overset{\tau}{\longrightarrow}
H^{0}(S;\OS(aL+bC)(\nu))$ is injective when $c_{1}=0$ and
$\overset{3}{\oplus}H^{0}(\PS{3};\Op{3}(\nu-1 )) \overset{\tau}{\longrightarrow}
H^{0}(S;\OS(aL+bC)(\nu-1))$ is injective when $c_{1}=-1,-2.$ For
definiteness, consider the first case, and note that it is
equivalent to showing that
$\overset{3}{\oplus}H^{0}(S;\OS(\nu )) \overset{\tau}{\longrightarrow}
H^{0}(S;\OS(aL+bC)(\nu))$ is injective. Set $V\equiv {\cnums}^{3}$
and extend the sheaf homomorphism
$V\otimes \OS(\nu) \overset{\tau}{\longrightarrow}\OS(aL+bC)(\nu)$
to a Koszul sequence over S

\begin{multline} \label{E:eq5}
0 \longrightarrow \wedge^{3}V\otimes \OS(-2aL-2bC)(\nu)
\longrightarrow \wedge^{2}V\otimes \OS(-aL-bC)(\nu) \\
\longrightarrow V\otimes \OS(\nu) \longrightarrow \OS(aL+bC)(\nu)
\longrightarrow 0.
\end{multline}

\noindent Break this up into two short exact sequences

\begin{equation} \label{E:eq6}
0 \longrightarrow \wedge^{3}V\otimes \OS(-2aL-2bC)(\nu)
\longrightarrow  \wedge^{2}V\otimes \OS(-aL-bC)(\nu)\longrightarrow \mK \longrightarrow 0
\end{equation}

\begin{equation} \label{E:eq7}
0 \longrightarrow \mK \longrightarrow V\otimes \OS(\nu) \longrightarrow
\OS(aL+bC)(\nu) \longrightarrow 0.
\end{equation}

\noindent Considering the second sequence, we must show that
$H^{0}(S;\mK)=0.$ By the first sequence and the fact that
$\OS(-aL-bC)(\nu)\cong \OS(-(a-b)L)(\nu-b)$ has no global
sections (Lemma~\ref{Lem:Main Lemma}i), it is enough to prove
that $H^{1}(S;\OS(-2aL-2bC)(\nu))\cong H^{1}(S;\OS(-2(a-b)L)(\nu
-2b))=0.$ Except in the cases $k=3$, $a=2l$, $b=l$, for $l\geq 2$,
this follows from $\nu-2b<0$ and
Lemma~\ref{Lem:Main Lemma}iv since $2b-\nu >(k-2)(2(a-b)-1)$
reduces to $2[b(k-1)-a(k-2)]+2\nu -2>0$ which holds since
$b(k-1)-a(k-2)\geq 0$ and $2\nu -2\geq 0$ and
both equalities hold iff $\nu=1$, k=3, and $a=2l$, $b=l$ for $l\in
{\znums}_{+}.$ When $k=3$, $a=2l$, $b=l$ for $l\geq 2$,  $h^{1}(S;\OS(-2aL-2bC)(\nu))
=h^{1}(S;\OS(-4lL-2lC)(1))=1$ but $H^{1}(S;\OS(-4lL-2lC)(1))
\rightarrow \overset{3}{\oplus}H^{1}(S;\OS(-2lL-lC)(1))$ is injective so that we
again get $H^{0}(S;\mK)=0.$ To see this, note that
$\OS(-4lL-2lC)(1)\cong \OS(-2lL)(1-2l)$ and $\OS(-2lL-lC)(1)\cong \OS(-lL)(1-l)$
and consider

$$
\begin{CD}
0 @>>> \overset{3}{\oplus}\OS(-lL)(1-l) @>>> \overset{3}{\oplus}\OS(-(l-1)L)(1-l)
@>>> \overset{3}{\oplus}\mO_{L} @>>> 0 \\
@AAA  @A\tau AA  @A\tau AA  @A\psi AA  @AAA  \\
0 @>>> \OS(-2lL)(1-2l) @>>> \OS(-(2l-1)L)(1-2l)
@>>> \mO_{L} @>>> 0 \\
\end{CD}
$$

\noindent and note that the homomorphism $\psi$ induced by $\tau$ is given
by a non-zero $v\in {\cnums}^{3}.$ The cohomology ladder and Lemma~\ref{Lem:Main Lemma},
for $l\geq 2$ gives

$$
\begin{CD}
\overset{3}{\oplus}H^{0}(L;\mO_{L}) \cong {\cnums}^{3} @>\cong >>
\overset{3}{\oplus}H^{1}(S;\OS(-lL)(1-l)) \\
@AvAA @A\tau AA \\
H^{0}(L;\mO_{L}) \cong \cnums @>\cong >> H^{1}(S;\OS(-2lL)(1-2l))
\end{CD}
$$

\noindent and this yields our result. The $c_{1}=-1,-2$ case follows in the same way.

Now assume $k\geq 3$ and $b\geq a$ (case 4.). Then (\ref{E:eq3})
requires $b>2\nu +c_{1}$ and (\ref{E:eq2}) implies $a\geq 0.$ Note
that $\OS(aL+bC)\cong \OS((b-a)C)(a) \cong \OS(a)\otimes
\OS((b-a)C)$ is globally generated because both $\OS(a)$ and $\OS((b-a)C)$
are (Lemma~\ref{Lem:Main Lemma}iii). To verify that $a>\nu /2$
(when $c_{1}=0$) implies that
$\overset{3}{\oplus}H^{0}(\PS{3};\Op{3}(\nu )) \overset{\tau}{\longrightarrow}
H^{0}(S;\OS(aL+bC)(\nu))$ is injective and that $a>(\nu-1)/2$
(when $c_{1}=-1, -2)$ implies
$\overset{3}{\oplus}H^{0}(\PS{3};\Op{3}(\nu-1)) \overset{\tau}{\longrightarrow}
H^{0}(S;\OS(aL+bC)(\nu-1))$ is injective, proceed as in case 3.
In the $c_{1}=0$ situation, this reduces to knowing that $H^{1}(S;\OS(-2(b-a)C)(\nu
-2a))=0.$ By Lemma~\ref{Lem:Main Lemma}v this holds if $a>\nu /2$
or $a=\nu /2=b.$ The second case cannot occur since $b>2\nu \geq \nu /2.$
The $c_{1}=-1,-2$ case is similar.

In case 2., $S=Q\cong \PS{1}\times \PS{1}$ , the argument follows the
same pattern and is left to the reader.
\end{proof}

\begin{Prop} \label{Prop:properties of rank 3 bundles}
Let E be a rank 3 stable bundle $on \PS{3}$ of the
type constructed in Theorem~\ref{Thm:main theorem rank 3} for $k\geq 3.$
Then
\begin{enumerate}
\item For $c_{1}=0$, $l\geq -4$, $H^{3}(\PS{3};E(l))=0.$ For $c_{1}=-1,-2$
, $l\geq -3$, $H^{3}(\PS{3};E(l))=0.$
\item For $l>\nu-4$, $H^{2}(\PS{3};E(l))=0$ iff $l>min(a,b)+\nu -4.$
\item $H^{1}(\PS{3};E(l))=0$ in exactly the following cases:For $b=a.$
For $b>a$, $l>b+\nu -4 +(k-2)[b-a-1].$ For $b=a+1$, $l=b+\nu -4.$ For $a>b$
, $l>a+\nu -4.$ For $l=a+\nu -4$, $a=b+1.$
\item For $a\geq b$ , E(l) is globally generated iff $l\geq max(a-k+\nu,\nu).$
For $b>a$, E(l) is globally generated iff $l\geq b(k-1)-a(k-2)-k+\nu.$
\end{enumerate}
\end{Prop}

\begin{proof}
This is very similar to the proof of Proposition~\ref{Prop:properties of
rank 2 bundles} and so is left to the reader.
\end{proof}

\bigskip

To count moduli, we proceed as in the previous examples. Fix
$\nu$, $c_{1}$, a, and b and consider
the dual defining sequences for E

\begin{equation} \label{E:ds1}
\begin{CD}
0 @>>> \overset{3}{\oplus}\Op{3}(-\nu) @>> \sigma > E
@>>> {j_{S}}_{*}\OS (-aL-bC)(2\nu+c_{1}) @>>> 0
\end{CD}
\end{equation}

\begin{equation} \label{E:ds2}
\begin{CD}
0 @>>> E^{*} @>> \sigma^{t} > \overset{3}{\oplus}\Op{3}(\nu) @>> \tau >
{j_{S}}_{*}\OS (aL+bC)(\nu) @>>> 0.
\end{CD}
\end{equation}

\noindent Here $\sigma \in \overset{3}{\oplus}H^{0}(\PS{3};E(\nu))$
and $\tau \in \overset{3}{\oplus}H^{0}(S;\OS(aL+bC))$.
The sequences imply that the function $(S,L,\tau)\rightarrow (E,\sigma)$
is injective but not a priori surjective (when $b\leq k$), as explained in
Section~\ref{S:examples of rank 2}: For E of the form
(\ref{E:ds1}), choose a different
$\bar{\sigma}\in \overset{3}{\oplus}H^{0}(\PS{3};E(\nu))$; this
gives another sequence

\begin{equation} \label{E:ds3}
\begin{CD}
0 @>>> \overset{3}{\oplus}\Op{3}(-\nu) @>> \bar{\sigma} > E
@>>> {j_{\bar{S}}}_{*}\bar{\mL} @>>> 0
\end{CD}
\end{equation}

\noindent and we must show that the surface $\bar{S}$ contains a
line $\bar{L}$ and that $\bar{\mL}$ has the form
$\mO_{\bar{S}} (-a\bar{L}-b\bar{C})(2\nu+c_{1}).$ This is obvious
when S and $\bar{S}$ are quadrics ($k=2$) so assume $k\geq 3.$

Consider case 4 of Theorem~\ref{Thm:main theorem rank 3}: $b\geq
a$, $a\geq \nu/2$ (for $c_{1}=0$), $\geq (\nu-1)$
(for $c_{1}=-1,-2$), $b>2\nu+c_{1}.$ Arguing exactly as in
Section~\ref{S:examples of rank 2}, one finds
$\bar{\mL}(b-2\nu-c_{1})\cong \mO_{\bar{S}}(\bar{D})$
for $\bar{D}$ effective, $deg\bar{D}=b-a$,
${\bar{D}}^{2}=-(k-2){(b-a)}^{2}$, and finally that $\bar{D}=(b-a)\bar{L}$
for $\bar{L}\subset \bar{S}$ a line.

Now consider case 3 of Theorem~\ref{Thm:main theorem rank 3}:
$a>b$, $b(k-1)-a(k-2)\geq 0.$ Using $b\leq k$, there are only two
possibilities: $a=k,b=k-1$ and $a=k-1,b=k-2.$ Again proceeding as
in Section~\ref{S:examples of rank 2}, $deg\bar{D}=k-1$ and
${\bar{D}}^{2}=0.$

When $a=k,b=k-1$, the cohomology sequences of (\ref{E:ds1}) and
(\ref{E:ds3}) and Lemma~\ref{Lem:Main Lemma} imply
$h^{0}(\bar{S};\bar{\mL}(\nu))=h^{0}(S;\OS(C))=2.$ Therefore
$\bar{\mL}(\nu)\cong \mO_{\bar{S}}(\bar{D})$ for $\bar{D}$ effective.
Since $\OS(C)$ is globally generated, so is $E(\nu)$ and therefore
$\mO_{\bar{S}}(\bar{D}).$
So $|\bar{D}|$ has no base locus and, by Bertini, we can assume $\bar{D}$
is smooth. $\bar{D}=\sum_{i}Y_{i}$ for the $Y_{i}$ smooth disjoint
curves and ${Y_{i}}^{2}=0.$ The cohomology sequence of

$$
\sexs{\mO_{\bar{S}}}{\mO_{\bar{S}}(\bar{D})}{\bigoplus_{i=1}^{r}\mO_{Y_{i}}}
$$

\noindent and $h^{0}(\bar{S};\mO_{\bar{S}}(\bar{D}))=2$ imply $r=1.$
So $\bar{D}$ is an irreducible smooth curve of degree $k-1.$ To
show that $\bar{D}$ is contained in a plane, we show that,
for $\Gamma_{1}$ denoting the homogeneous polynomials in x of degree 1, the
restriction $\Gamma_{1}\rightarrow H^{0}(\bar{D};\mO_{\bar{D}}(1))$
has non-trivial kernel. This follows from the cohomology sequence
of

$$
\sexs{\mO_{\bar{S}}(-\bar{D})(1)}{\mO_{\bar{S}}(1)}{\mO_{\bar{D}}(1)}
$$

\noindent and $h^{0}(\bar{S};\mO_{\bar{S}}(-\bar{D})(1))=h^{0}(S;\OS(-C)(1))=
h^{0}(S;\OS(L))=1.$ So $\bar{S}$ has a hyperplane section $H=\bar{D}+\bar{L}$
for $\bar{L}$ a line. This gives
$\bar{\mL} \cong \mO_{\bar{S}}(-k\bar{L}-(k-1)\bar{D})(2\nu+c_{1}).$

When $a=k-1,b=k-2$, the argument is the same except that we must
work harder to show that $|\bar{D}|$ has no base locus. The
sequences (\ref{E:ds1}) and (\ref{E:ds3}) show that
$\bar{\mL}(\nu-1)\cong \mO_{\bar{S}}(\bar{D})$ for $\bar{D}$
effective and that $h^{0}(\bar{S};\mO_{\bar{S}}(\bar{D}))=2.$ Note that
$\mO_{\bar{S}}(\bar{D})(1)$ is globally generated but that $\mO_{\bar{S}}(\bar{D})$ is
not. Therefore to prove that $|\bar{D}|$ has no base locus it is enough to
show that the bilinear multiplication map

$$
m:\Gamma_{1}\times H^{0}(\bar{S};\mO_{\bar{S}}(\bar{D}))\rightarrow
H^{0}(\bar{S};\mO_{\bar{S}}(\bar{D})(1))
$$

\noindent is surjective. Since
$h^{0}(\bar{S};\mO_{\bar{S}}(\bar{D})(1))=7$, we must show that
the dimension of the kernel of m is $\leq 1.$ For $V\cong {\cnums}^{4}$
defined by $\PS{3}\equiv \mP V$ and x the homogeneous coordinates
on $\PS{3}$, consider

\begin{equation} \label{E:ds4}
\begin{CD}
0 @>>> \mK_{1} @>>> \mO_{\bar{S}}(\bar{D})\otimes V @>> \iota_{x} >
\mO_{\bar{S}}(\bar{D})(1) @>>> 0.
\end{CD}
\end{equation}

\noindent and note that $\Ker m \cong H^{0}(\bar{S};\mK_{1}).$
Extend (\ref{E:ds4}) to a Koszul sequence which breaks up into
three short exact sequences, (\ref{E:ds4}) and

$$
\sexs{\mK_{2}}{\mO_{\bar{S}}(\bar{D})(-1)\otimes
{\wedge}^{2}V}{\mK_{1}}
$$

$$
\sexs{\mO_{\bar{S}}(\bar{D})(-3)\otimes {\wedge}^{4}V}
{\mO_{\bar{S}}(\bar{D})(-2)\otimes {\wedge}^{3}V}{\mK_{2}}.
$$

\medskip

\noindent $h^{i}(\bar{S};\mO_{\bar{S}}(\bar{D})(-j))=
h^{i}(S;\OS(C)(-j))=h^{i}(S;\OS(-L)(-(j-1)))=0$ for $i=0,1$ and $j=1,2,3.$
by Lemma~\ref{Lem:Main Lemma}iv. Therefore dim $\Ker m$ equals
the dimension of the kernel of

$$
H^{2}(\bar{S};\mO_{\bar{S}}(\bar{D})(-3))\otimes
{\wedge}^{4}V \overset{\iota_{x}}{\longrightarrow}
H^{2}(\bar{S};\mO_{\bar{S}}(\bar{D})(-2))\otimes
{\wedge}^{3}V
$$

\noindent which, by Serre duality and (\ref{E:ds1}),
(\ref{E:ds3}), equals the dimension of the cokernel of

\begin{equation}  \label{E:ds5}
H^{0}(S;\OS(-C)(k-2))\otimes {\wedge}^{3}V \overset{x\wedge}{\longrightarrow}
H^{0}(S;\OS(-C)(k-1))\otimes {\wedge}^{4}V.
\end{equation}

\noindent Using $\OS(-C)(k-2)\cong \OS(L)(k-3)$ and the sequence

$$
\sexs{\OS(k-3)}{\OS(L)(k-3)}{\mO_{L}(-1)}
$$

\noindent we see that $\OS(-C)(k-2)$ is not globally generated because all
its sections vanish on L. Similarly, $\OS(-C)(k-1)\cong \OS(L)(k-2)$
is globally generated and we can take a basis for $H^{0}(S;\OS(L)(k-2))$
of the form $s_{0}, s_{1}, \dots, s_{n}$ where $s_{0}$ is
non-vanishing on L and the $s_{i}$ are zero on L for i=1 to n. Let
$\xi=0$ define L. Then the sequence shows that $s\in H^{0}(S;\OS(L)(k-3))$
has the form $s=\xi P$ for $P\in \Gamma_{k-3}.$ Similarly,
$\tilde{s}\in span\{s_{1}, \dots, s_{n}\}$ has the form $\tilde{s}=\xi \tilde{P}$
for $\tilde{P}\in \Gamma(k-2).$ This shows that the cokernel of (\ref{E:ds5})
has dimension 1.

We have demonstrated the 1-to-1 correspondences

$$
(E,\sigma)\longleftrightarrow (S,L,\tau)
$$

\noindent for $k\geq 3$ and

$$
(E,\sigma) \longleftrightarrow (Q,\pm,\tau) \text{\quad for } k=2.
$$

\begin{Prop}  \label{Prop:parameters rank 3}
Let $\mY$ be the set of isomorphism classes of stable rank 3
bundles E as in Theorem~\ref{Thm:main theorem rank 3}.

\noindent For $b\geq a\geq k-3$,

\begin{align}
dim\mY&=3(k-1)ab-\frac{3(k-2)}{2}a^{2}-\frac{3(k-4)}{2}(a+(k-1)b)
\notag \\
\quad &+3\binom{k-1}{3}+\binom{k+3}{3}-sup(k-3,0)-7 \notag \\
\quad &+(\text{when }b\leq k)\begin{cases}
-3\binom{k-b+3}{3} \quad&\text{if }k\geq 3 \\
3a-9 \quad &\text{if }k=2.  \notag
\end{cases}
\end{align}

\noindent For $a>b$,

\begin{align}
dim\mY&=3(k-1)ab-\frac{3(k-2)}{2}a^{2}-\frac{3(k-4)}{2}(a+(k-1)b)
\notag \\
\quad &+3\binom{k-1}{3}+\binom{k+3}{3}-sup(k-3,0)-7 \notag \\
\quad &+(\text{when }a\leq k)\begin{cases}
-6 \quad&\text{if }a=k \\
-21\quad &\text{if }a=k-1.  \notag
\end{cases}
\end{align}

\noindent For $b\geq a$, $a\leq k-2$,

\begin{align}
dim\mY &=\binom{a+2}{2}[3b-2a+3]+\binom{k+3}{3}-sup(k-3,0)-10
\notag \\
\quad &+(\text{when } b\leq k) \begin{cases}
-3\binom{k-b+3}{3} &\text{if }k\geq 3 \\
-9 &\text{if } k=2.  \notag
\end{cases}
\end{align}

\end{Prop}

\begin{proof}
Counting parameters from the 1-to-1 correspondences
and using $dim\{(S,L)\}=dim\{S\}$ when $k\geq
3$, we get

\begin{align} \label{E:e7}
dim \mY &=dim\{S\}+dim\{\tau\}-dim\{\sigma\}     \\
&=\binom{k+3}{3}-\sup(k-3,0)-10+3h^{0}(S;\OS(aL+bC))  \notag \\
& \quad -3h^{0}(S;\OS(-aL-bC)(k)). \notag
\end{align}

\noindent When $max(a,b)>k$, $h^{0}(S;\OS(-aL-bC)(k))=0$ and $h^{0}(S;\OS(aL+bC))$
is given by Lemma~\ref{Lem:sections}. When $max(a,b)\leq k$,
Lemma~\ref{Lem:sections} can also
be applied to
$h^{0}(S;\OS(-aL-bC)(k))=h^{0}(S;\OS((k-a)L-(k-b)C))$, using the
restrictions on a and b imposed by Theorem~\ref{Thm:main theorem
rank 3}. This gives our result.
\end{proof}

\bigskip

From

\begin{equation} \label{E:s3}
\begin{CD}
0 @>>> \overset{3}{\oplus}\Op{3} @>> \sigma > E(\nu)
@>>> {j_{S}}_{*}\OS (-aL-bC)(k) @>>> 0
\end{CD}
\end{equation}

\noindent it follows that

\begin{equation} \label{E:e4}
H^{j}(\PS{3};E(\nu)) \cong H^{j}(S;\OS (-aL-bC)(k)) \quad j=1,2,3
\end{equation}

\begin{align} \label{E:e5}
h^{0}(\PS{3};E(\nu)) &= 3+ h^{0}(S;\OS (-aL-bC)(k))  \\
\quad &=3 \quad \text{if } max(a,b)>k \notag \\
\quad &\text{(by Lemma~\ref{Lem:Main Lemma}i and ii).}  \notag
\end{align}

\bigskip

To establish a framework in which we can try to
calculate or estimate the dimension of the Zariski tangent
space of $\mM$ at E we argue as in Section~\ref{S:examples of rank
2} to obtain

\begin{equation} \label{E:e8}
\sexs{{\sheafend}(E)}{\overset{3}{\oplus}E(\nu)}{{j_{S}}_{*}E_{S} (aL+bC)(\nu)}
\end{equation}

\begin{multline} \label{E:e9}
0\longrightarrow \OS \longrightarrow \overset{3}{\oplus}\OS(aL+bC)
\longrightarrow {\mathcal K} \longrightarrow 0  \\
0 \longrightarrow {\mathcal K} \longrightarrow E_{S}(aL+bC)(\nu ) \longrightarrow \OS(k)
\longrightarrow 0
\end{multline}

\begin{multline} \label{E:e10}
0 \longrightarrow H^{0}(\PS{3};\sheafend(E)) \longrightarrow
\overset{3}{\oplus}H^{0}(\PS{3};E(\nu )) \longrightarrow
H^{0}(S;E_{S}(aL+bC)(\nu )) \\
\longrightarrow H^{1}(\PS{3};\sheafend(E))
\longrightarrow \overset{3}{\oplus}H^{1}(\PS{3};E(\nu ))
\longrightarrow H^{1}(S;E_{S}(aL+bC)(\nu )) \\
\longrightarrow H^{2}(\PS{3};\sheafend(E)) \longrightarrow
\overset{3}{\oplus}H^{2}(\PS{3};E(\nu )) \longrightarrow H^{2}(S;E_{S}(aL+bC)(\nu ))
\longrightarrow 0.
\end{multline}

The main result of the following calculations will be to make an
effective comparison of $dim T\mM_{E} =h^{1}(\PS{3};\sheafend E)$ and $dim \mY$
when k=2 or 3 and an estimation of the codimension of $\mY$ in $\mM$
at E when $k\geq 4.$
Using Lemma~\ref{Lem:Main Lemma} we get

\begin{align} \label{E:e11}
h^{1}(S;\OS(aL+bC)) &=0 \quad \text{in case 3.}  \\
\quad &=0 \quad \text{in case 4. iff } a>k-4 \text{ or }
a=k-4, b=a, a+1  \notag \\
h^{2}(S;\OS(aL+bC)) &=0 \quad \text{iff } max(a,b)>k-4.
\end{align}

\noindent From the cohomology sequences of (\ref{E:e9}) it follows that

\begin{align} \label{E:e13}
h^{2}(S;E_{S}(aL+bC)(\nu))&=0 \quad \text{for } max(a,b)> k-4   \\
h^{1}(S;E_{S}(aL+bC)(\nu))&= dim(coker\delta)  \notag \\
\quad &=0 \text{ for } k=2,3 \notag \\
h^{0}(S;E_{S}(aL+bC)(\nu))&=\binom{k+3}{3}-2+3h^{0}(S;\OS(aL+bC))-dim(im\delta)
 \notag
\end{align}

\noindent where $\delta$ is the connecting homomorphism
$H^{0}(S;\OS(k))\rightarrow H^{1}(S;\mK).$ When $k=2,3$ ,  $\delta =0$
because $H^{1}(S;\mK)=0.$
Starting from (\ref{E:e4}), routine calculations using
Lemma~\ref{Lem:Main Lemma} give

\begin{align} \label{E:e14}
h^{2}(\PS{3};E(\nu))&=h^{2}(S;\OS (-aL-bC)(k)) \\
\quad &=h^{0}(S;\OS (aL+bC)(-4)) \notag \\
\quad &=h^{0}(S;\OS ((a-4)L+(b-4)C))  \text{ (see Lemma~\ref{Lem:sections} when }
min(a,b)\geq 4) \notag \\
\quad &=0 \text{ if } min(a,b)<4  \notag
\end{align}

\noindent and

\begin{align} \label{E:e19}
h^{1}(\PS{3};E(\nu))&=h^{1}(S;\OS (-aL-bC)(k))  \\
&=0 \text{ in case 3 of Theorem~\ref{Thm:main theorem rank 3}
if } b\geq k \text{ and } (k-1)b-(k-2)a>2 \text{ or }  \notag  \\
&a=k+1,b=k \text{ and in case 4 of Theorem~\ref{Thm:main theorem rank 3}
if } a>k \notag  \\
&\text{ or } a=k, b=k,k+1.   \notag
\end{align}

\bigskip

\begin{Thm} \label{Thm:moduli for k large}
Let$S\hookrightarrow \PS{3}$ be a smooth surface of degree k that contains a line L and
let E be a stable rank 3 bundle on $\PS{3}$ of the form

\begin{equation}
\sexs{\overset{3}{\oplus}\Op{3}(-\nu)}{E}{{j_{S}}_{*}\OS(-aL-bC)(2\nu +c_{1})}
\notag
\end{equation}

\noindent as in Theorem~\ref{Thm:main theorem rank 3}. Assume in
case 3 that $b\geq k$ and $b(k-1)-a(k-2)>2$ or that $a=k+1,b=k$ and in case 4 that $a>k$
or $a=k$, $b=k+1.$ Then

\begin{equation} \label{E:e20}
h^{1}(\PS{3};\sheafend E)=\binom{k+3}{3}-10+3h^{0}(S;\OS(aL+bC))-dim(im\delta).
\end{equation}

\noindent For $\mY \subset \mM$ the set of these bundles, $\mY$ is a subscheme of
$\mM$ of codimension no larger than
$max(k-3,0)-dim(im\delta)\leq max(k-3,0).$
\end{Thm}

\begin{proof}
By (\ref{E:e19}), our hypothesis insures that $H^{1}(\PS{3};E(\nu))=0$
and that \linebreak
$H^{0}(S;\OS(-aL-bC)(k))=0.$ Now the sequence (\ref{E:e10}) gives

$$
h^{1}(\PS{3};\sheafend E)=h^{0}(S;E_{S}(aL+bC)(\nu
))-3h^{0}(\PS{3};E(\nu ))+1  \notag
$$

\noindent and the formula follows from (\ref{E:e13}) and
(\ref{E:e5}). Comparing this with dim $\mY$ (given by
(\ref{E:e7}))and using $dim_{E}\mM \leq h^{1}(\PS{3};\sheafend E)$
gives the codimension estimate. $\mY$ is a subscheme
of $\mM$ from arguments parallel to those in
Section~\ref{S:scheme}.
\end{proof}

\bigskip

We now assume k=2 or 3 and obtain much more precise information.
It will be convenient to calculate $h^{1}(\PS{3};\sheafend E)$ via
$h^{1}=[h^{1}-h^{2}]+h^{2}.$ From Riemann-Roch (\ref{E:exp dim of
moduli rank 3})and the chern class formula of
Theorem~\ref{Thm:main theorem rank 3},

\begin{align} \label{E:e21}
h^{1}(\PS{3};\sheafend E)-h^{2}(\PS{3};\sheafend E)&=12c_{2}(E)-4c_{1}^{2}-8
\notag \\
\quad &=12[a+b(k-1)]-4k^{2}-8.
\end{align}

\noindent From (\ref{E:e10}), (\ref{E:e13}), and (\ref{E:e14}),

\begin{align} \label{E:e22}
h^{2}(\PS{3};\sheafend E)&=3h^{0}(S;\OS((a-4)L+(b-4)C))
\text{ (see Lemma~\ref{Lem:sections} when }
min(a,b)\geq 4) \notag \\
\quad &=0 \text{ if } min(a,b)<4   \notag \\
\quad &\text{ for k=2 or 3.}
\end{align}

\noindent Our conclusions are summarized in the following two
theorems.

\bigskip

\begin{Thm} \label{Thm:moduli degS=2}
Let $Q\hookrightarrow \PS{3}$ be a smooth quadric and let
E be a stable rank 3 bundle on $\PS{3}$ of the form

\begin{equation}
\sexs{\overset{3}{\oplus}\Op{3}(-1)}{E}{{j_{Q}}_{*}\mO_{Q}(1-a,1-b)}
\notag
\end{equation}

\noindent with $a,b\geq 0, max(a,b)\geq 2$ as in Theorem~\ref{Thm:main theorem rank
3}. Let $\mY \subset \mM$ be the set of these bundles.Then

\begin{equation}
dim \mY=
\begin{cases}
3(a+1)(b+1) &\quad \text{for } max(a,b)\geq 3  \\
12a+12b-24 &\quad \text{for } max(a,b)\leq 3   \notag
\end{cases}
\end{equation}

\begin{equation}
h^{1}(\PS{3};\sheafend E)=
\begin{cases}
3(a+1)(b+1) &\quad \text{for } min(a,b)\geq 3  \\
12a+12b-24 &\quad \text{for } min(a,b)\leq 3   \notag
\end{cases}
\end{equation}

\begin{equation}
h^{2}(\PS{3};\sheafend E)=
\begin{cases}
3(a-3)(b-3) &\quad \text{for } min(a,b)\geq 3 \\
0 &\quad \text{for } min(a,b)\leq 3.    \notag
\end{cases}
\end{equation}

\noindent In all cases $dim_{E}\mM =h^{1}(\PS{3};\sheafend E)= dim T\mM_{E}$
and $\mM$ is smooth at E. $dim_{E}\mM =dim \mY$ in all cases except when
$max(a,b)\geq 4$ and $min(a,b)\leq 2$ both hold.
\end{Thm}

\begin{proof}
For $k=2$, Lemma~\ref{Lem:sections} reduces to the elementary result

$$
h^{0}(Q;\mO_{Q}(i,j))=
\begin{cases}
(i+1)(j+1) &\text{ if } i,j\geq 0 \notag \\
0 &\text{ if } i \text{ or } j<0.
\end{cases}
$$

\noindent This along with (\ref{E:e7}), (\ref{E:e21}), and
(\ref{E:e22}) give our three formulas.
In every case, either $h^{1}(\PS{3};\sheafend E)=dim\mY$
or $h^{2}(\PS{3};\sheafend E)=0.$
\end{proof}

\bigskip

\begin{Thm} \label{Thm:moduli degS=3}
Let $S_{3}\hookrightarrow \PS{3}$ be a smooth cubic and let E
be a stable rank 3 bundle on $\PS{3}$ of the form

\begin{equation}
\sexs{\overset{3}{\oplus}\Op{3}(-1)}{E}{{j_{S_{3}}}_{*}\mO_{S_{3}}(-aL-bC)(2)}
\notag
\end{equation}

\noindent as in Theorem~\ref{Thm:main theorem rank 3}. Let $\mY \subset \mM$
be the set of these bundles. Then
\begin{equation}
dim \mY=
\begin{cases}
6ab+3b-3\binom{a}{2}+13 &\quad \text{for } max(a,b)\geq 4  \\
12a+24b-44 &\quad \text{for } max(a,b)\leq 3,a\geq b   \notag \\
52 &\quad b=3,a=2 \notag \\
37 &\quad b=3,a= 1 \notag
\end{cases}
\end{equation}

\begin{equation}
h^{2}(\PS{3};\sheafend E)=
\begin{cases}
0 &\quad \text{for } min(a,b)\leq 3  \notag \\
3[2ab-\frac{1}{2}a^{2}-\frac{7}{2}a-7b+19] &\quad \text{if }
min(a,b)\geq 4 \text{ and } \\
\quad &\quad b\geq a \text{ or } a>b\geq 2+a/2 \text{ both hold } \\
3[2b^{2}-14b+25] &\quad \text{if } min(a,b)\geq 4 \text{ and } \\
\quad &\quad b<2+a/2 \text{ both hold }
\end{cases}
\end{equation}

\begin{equation}
h^{1}(\PS{3};\sheafend E)=
\begin{cases}
12a+24b-44 &\quad \text{for } min(a,b)\leq 3  \notag \\
6ab+3b-3\binom{a}{2}+13 &\quad \text{if } min(a,b)\geq 4  \text{ and }
\\
\quad &\quad b\geq a \text{ or } a>b\geq 2+a/2 \text{ both hold }\\
12a+6b^{2}-18b+31 &\quad \text{if } min(a,b)\geq 4 \text{ and } \\
\quad &\quad b<2+a/2 \text{ both hold. }
\end{cases}
\end{equation}

\noindent [Note that in the case where $min(a,b)\geq 4$ and $b<2+a/2$
both hold, the requirement $2b-a\geq 0$ from Theorem~\ref{Thm:main
theorem rank 3} gives that $b=a/2$ or $a/2+1$ when a is even and
$b=(a+1)/2$ or $(a+1)/2+1$ when a is odd.]

\noindent $dim\mY=dim_{E}\mM=h^{1}(\PS{3};\sheafend
E)$ and $\mM$ is smooth at E
when $min(a,b)\geq 4$ and $b\geq a$ or $a>b\geq 2+a/2$ both hold,
when $max(a,b)\leq 3, a\geq b$, when $b=3,a=2$, and when $min(a,b)\geq 4, b=a/2+3/2$
for a odd. When $min(a,b)\leq 3$, $\mM$ is smooth at E.
\end{Thm}

\begin{proof}
This follows from
Lemma~\ref{Lem:sections}, Proposition~\ref{Prop:parameters rank 3}, (\ref{E:e21}), and
(\ref{E:e22}).
\end{proof}

\bigskip

\nocite{Maruyama2}
\nocite{Maruyama3}
\nocite{Gieseker-ms}

\providecommand{\bysame}{\leavevmode\hbox to3em{\hrulefill}\thinspace}
\providecommand{\MR}{\relax\ifhmode\unskip\space\fi MR }
\providecommand{\MRhref}[2]{%
  \href{http://www.ams.org/mathscinet-getitem?mr=#1}{#2}
}
\providecommand{\href}[2]{#2}

\end{document}